\documentclass{article}

\usepackage{arxiv}

\usepackage{geometry}
\usepackage{graphicx}%
\usepackage{multirow}%
\usepackage{amsmath,amssymb,amsfonts}%
\usepackage{amsthm}%
\usepackage{mathrsfs}%
\usepackage[title]{appendix}%
\usepackage{xcolor}%
\usepackage{textcomp}%
\usepackage{manyfoot}%
\usepackage{booktabs}%
\usepackage{algorithm}%
\usepackage{algorithmicx}%
\usepackage{algpseudocode}%
\usepackage{listings}%
\allowdisplaybreaks[4]

\theoremstyle{thmstyleone}%
\newtheorem{theorem}{Theorem}
\theoremstyle{thmstyletwo}%
\newtheorem{lemma}{Lemma}
\newtheorem{expl}{Example}

\theoremstyle{thmstylethree}%

\newcommand{\EE}{\mathbb{E}}
 \newcommand{\PP}{\mathbb{P}}
 \newcommand{\RR}{\mathbb{R}}

 \newcommand{\ind}{\mathbf{1}}

 \newcommand{\eproof}{\indent\vrule height6pt width4pt depth1pt\hfil\par\medbreak}

 \newtheorem{assumption}{Assumption}

 \def\F{{\cal F}}

 \def\a{\alpha}
 \def\m{\mu}
 \def\k{\kappa}
 \def\r{\rho}

\title{A Milstein-type method for highly non-linear non-autonomous time-changed stochastic differential equations}

\author{
 Wei Liu, Ruoxue Wu, Ruchun Zuo \\
 Department of Mathematics, Shanghai Normal University, Shanghai, 200234, China\\
 \texttt{weiliu@shnu.edu.cn; 963729621@qq.com; zuoruchun@qq.com}
}

\begin{document}
\maketitle

\begin{abstract}
A Milstein-type method is proposed for some highly non-linear non-autonomous time-changed stochastic differential equations (SDEs). The spatial variables in the coefficients of the time-changed SDEs satisfy the super-linear growth condition and the temporal variables obey some H\"older's continuity condition. The strong convergence in the finite time is studied and the convergence order is obtained.
\end{abstract}


\section{Introduction}\label{sec1}

Time-changed stochastic processes and time-changed stochastic differential equations (SDEs) have been attracting increasing attentions in the past decades, as they are one of the important tools to describe sub-diffusion processes and their close relation with deterministic fractional differential equations (DFDE) \cite{UHK2018}. 

In \cite{MS2004}, Meerschaert and Scheffler gave a detailed discussion on the time process used for changing times and established a fundamental limit theorem that links some continuous-time random walks with infinite  mean waiting times with a class of time-changed stochastic processes. Some important properties and essential inequalities of time-changed fractional Brownian motion were obtained by Deng and Schilling in \cite{DS2017}. 

The existence and uniqueness theorem for time-changed SDEs and many useful tools were obtained by Kobayashi in \cite{Kob2011}.  Stabilities in different senses  of all kinds of stochastic equations were broadly discussed: Wu in \cite{Wu2016} investigated SDEs driven time-changed Brownian motion; Nane and Li in \cite{NN2017,NN2018} studied the case when the driven noise is the time-changed  L\'evy noise;  Zhang and Yuan focused on time-changed stochastic functional differential equations in \cite{ZY2019}; Yin et al. considered the impulsive effects on stabilities in \cite{YXS2021} for a class of time-changed SDEs; Shen et al. in \cite{SZSW2023} discussed distribution dependent SDEs driven by time-changed Brownian motions.  Li et al. discussed some theoretical results of the time-changed McKean-Vlasov SDE in \cite{LXY2023}, which is also a distribution dependent SDE. 

Time-changed processes and time-changed SDEs are widely applied in modelling financial markets. Magdziarz introduced sub-diffusive Black-Scholes formula by using the classical geometric Brownian motion with the inverse $\alpha$-stable subordinator \cite{Mag2009a}.  Magdziarz et al. proposed the sub-diffusive version of the Bachelier Model and investigated its application in the option pricing \cite{MOW2011}. Janczura et al. studied the time-changed Ornstein-Uhlenbeck process that is driven by the $\alpha$-stable process and fitted the data from emerging markets in this model \cite{JOW2011}.  For connections between time-changed processes and various DFDE, we refer the readers to \cite{Che2017,HKU2012,Mag2009b,NN2016} and references therein.

Since the following two main reasons, numerical approximations to time-changed SDEs become essential. (1) Explicit forms of true solutions to time-changed SDEs are hardly found. (2) Applications of time-changed SDE models in practice often require a considerable number of sample paths to conduct statistical learnings like estimations, tests and predictions based on observed data. In this case, even explicit expressions of true solutions to some types of time-changed SDE  models are available, performing those calculations without the aid of computer simulations is highly unlikely.

When transition probabilities of solutions to time-changed SDEs are needed to be simulated, the typical approaches used are discretising the corresponding deterministic fractional differential equations. There are fruitful works on numerical methods for DFDE and a far-from-complete list of them includes \cite{DFF2002,DGLZ2012,LZ2015,WZ2023} among many others. 

\medskip 
\noindent
In this paper, we focus on another important aspect of numerical approximations  to time-changed SDEs, i.e. numerical simulations of sample paths of solutions. 

In this aspect, Kobayashi and collaborators studied different numerical methods for time-changed SDEs with different structures, when the global Lipschitz conditions are imposed on the spatial variables in the coefficients. The convergences in both the strong and weak senses of the Euler–Maruyama (EM) method for a class of time-changed SDEs were proved by Jum and Kobayashi in \cite{JK2016}, which, to our best knowledge, is the first work to study simulations of sample paths of solutions to time-changed SDEs. More recently, Jin and Kobayashi investigated some Euler-type and Milstein-type methods for more general type of time-changed SDEs in \cite{JK2019,JK2021}. One of the main differences in terms of techniques used  between \cite{JK2016} and \cite{JK2019,JK2021} is that the duality principle established in \cite{Kob2011} was employed in \cite{JK2016} but not in \cite{JK2019,JK2021}. Briefly speaking, the duality principle reveals the relation between the classical SDEs and time-changed counterpart, which enables numerical methods for time-changed SDEs to be constructed by using numerical methods for classical SDEs directly. For time-changed McKean-Vlasov SDEs,  Wen et al. considered the numerical method in \cite{WLX2023}.

In the case that some super-linear terms are allowed to appear in the coefficients, implicit methods and modified explicit methods are usually good alternatives as the classical Euler-type and Milstein-type methods may not be convergent \cite{HJK2011}. When some super-linear growth conditions are imposed on the spatial variables in the drift coefficient of time-changed SDEs, Deng and Liu studied the semi-implicit EM method in \cite{DL2020}, Liu et al. investigated the truncated EM method in \cite{LMTW2020} with the help of the duality principle, while Li et al. in \cite{LLLX2023} also discussed the truncated-type Euler method but without employing the duality principle.

In this paper, we also focus on numerical methods for time-changed SDEs with super-linear coefficients. Compared with \cite{DL2020,LMTW2020,LLLX2023}, we consider the numerical method with the higher convergence order by proposing a Milstein-type method with the truncating techniques to suppress super-linear terms. Due to the higher convergence order, compared with those Euler-type methods  Milstein-type methods are more suitable for the multi-level Monte Carlo that is quite popular for applications in finance  \cite{Gil2008,GS2013}.

\section{Mathematical preliminaries}

Let $(\Omega_W,\F^W,\PP_W)$ be a complete probability space with a filtration $\{\F^W_t\}_{t \geq 0}$ being right continuous and increasing, while $\F^W_0$ contains all $\PP_W$-null sets. Let $W(t)$ be a one-dimensional Wiener process defined in that probability space and is $\F^W_t$-adapted. $\EE_W$ denotes the expectation with respect to $\PP_W$. 

Let $(\Omega_D , \F^D, \PP_D)$ be another complete probability space  with a filtration $\left\{\F^D_t\right\}_{t \ge 0}$. $D(t)$ denotes a one-dimensional $\F^D_t$-adapted strictly increasing L\'evy process on $[0,\infty)$ starting from $D(0)=0$  defined on $(\Omega_D , \F^D, \PP_D)$.  Let $\EE_D$ denote the expectation with respect to $\PP_D$.  For detailed introductions and discussions on such a $D(t)$, we refer the readers to \cite{App2009,Sat1999} 

In this paper, $W(t)$ and $D(t)$ are assumed to be independent. Define the product probability space by $ (\Omega , \F, \PP):= (\Omega_W \times \Omega_D, \F^W\otimes \F^D, \PP_W \otimes \PP_D)$. Let $\EE$ denote the expectation under the probability measure $\PP$. It is clear that $\EE(\cdot) = \EE_D \left( \EE_W (\cdot) \right) = \EE_W  \left( \EE_D (\cdot) \right)$. 

For $x\in \RR^d$, $|x|$ denotes the Euclidean norm. The transposition of $x$ is denoted by $x^\mathrm{T}$. For two real numbers $a$ and $b$, set $a\vee b=\max(a,b)$ and $a\wedge b=\min(a,b)$. For a given set $G$, its indicator function is denoted by $\ind_{G}$.

Since $D(t)$ is strictly increasing, we define the inverse of $D(t)$ by
\begin{align*}
    E(t) := \inf\{ s\geq0\,;\,D(s) > t \}, \quad t \geq 0.
\end{align*}
Then, the $E(t)$ is used for changing time, as $t\mapsto E(t)$ is continuous and non-decreasing. The process $W(E(t))$ is called a time-changed Wiener process and $W(E(t))$ is regarded as a sub-diffusive process. For the simplicity of notations, we consider the one-dimensional $W(E(t))$ in our work. When $W(t)$ is a multi-dimensional Wiener process and the same $E(t)$ is used for changing time in each entry of $W(t)$, the results in this paper should still hold. But if different $E(t)$s are used to change times in different entries of $W(t)$, our results may not be applicable.

The time-changed SDEs considered in this paper take the following form, For any $T>0$ and $t\in [0,T]$
\begin{align}\label{SDE}
    dY(t) = f(t,Y(t))dE(t) + g(t,Y(t))dW(E(t)),~~~Y(0) = Y_0,
\end{align}
with $\EE |Y_0|^q < \infty$ for all $q > 0$, where $f: \RR_+ \times \RR^d \rightarrow \RR^d$ and $g: \RR_+ \times \RR^d \rightarrow \RR^{d}$.

Before we impose assumptions on the coefficients of \eqref{SDE}, we present some tedious but helpful notations. For any $y=(y^1,y^2,...,y^d) \in \RR^d$ and and any $t\in [0,T]$, define
\begin{align*}
    Lg(t,y)=\sum_{l=1}^{d}g^{l}(t,y)G^{l}(t,y),
\end{align*}
where $g=(g^{1},g^{2},...,g^{d})^{T}$, $g^{l}:\RR_+ \times \RR^d \rightarrow \RR$ and
\begin{align*}
    G^{l}(t,y)=\left( \dfrac{\partial g^{1}(t,y)}{\partial y^{l}},\dfrac{\partial g^{2}(t,y)}{\partial y^{l}},...,\dfrac{\partial g^{d}(t,y)}{\partial y^{l}} \right)^\mathrm{T}.
\end{align*}




The following assumptions are imposed on the coefficients of \eqref{SDE}. We first give requirements on spatial variables in the coefficients.
\begin{assumption}
    \label{ass1}
    Assume that there exist positive constants $\a$ and $C$ such that
    \begin{align*}
        |f(t,x)-f(t,y)|\vee |g(t,x)-g(t,y)|\vee |Lg(t,x)-Lg(t,y)|\leq C(1+|x|^{\a}+|y|^{\a})|x-y|,
    \end{align*}
for all $t\in[0,T]$ and any $x,y \in \RR^d$. 
\end{assumption}
It can be observed from Assumption \ref{ass1} that for all $t\in [0,T]$ and any $x\in \RR^d$
    \begin{align}
        \label{equ2-2}
        |f(t,x)|\vee|g(t,x)|\vee|Lg(t,x)|\leq M(1+|x|^{\a+1}),
    \end{align}
where $M$ depends on $C$ and $\sup_{0\leq t\leq T}\left(|f(t,0)|+|g(t,0)|+|Lg(t,0)|\right)$.

\begin{assumption}
    \label{ass2}
    Assume that there exists a pair of constants $p>2$ and $K>0$ such that
    \begin{align*}
        (x-y)^{\mathrm{T}}(f(t,x)-f(t,y)) +(5p-1)|g(t,x)-g(t,y)|^2\leq K|x-y|^2,
    \end{align*}
    for all $t\in [0,T]$ and any $x,y \in \mathbb{R}^d$.
\end{assumption}

\begin{assumption}
    \label{ass3}
    Assume that there exists a pair of constants $q>2$ and $K_1>0$ such that
    \begin{align*}
        x^{\mathrm{T}}f(t,x) +(5q-1)|g(t,x)|^2\leq K_1(1+|x|^2),
    \end{align*}
    for all $t\in [0,T]$ and any $x \in \mathbb{R}^d$.
\end{assumption}
Similar to the relation between Assumption \ref{ass1} and \eqref{equ2-2}, Assumption \ref{ass3} can be derived from Assumption \ref{ass2} but with complicated relations between $p$ and $q$ as well as $K$ and $K_1$. So we present Assumption \ref{ass3} as a new assumption.

\begin{assumption}
    \label{ass5}
    Assume that there exists a positive constant $ M^{\prime}$such that
    \begin{eqnarray*}
        &&|\frac{\partial f(t,x)}{\partial x}|\vee|\frac{\partial^{2} f(t,x)}{\partial x^{2}}|\vee|\frac{\partial g(t,x)}{\partial x}|\vee|\frac{\partial^{2} g(t,x)}{\partial x^{2}}|\leq M^{\prime}(1+|x|^{\a+1}),
    \end{eqnarray*}
    for any $x\in \mathbb{R}^d$ and all $t\in [0,T]$.
\end{assumption}

Now we turn to the requirement on the temporal variables in the coefficients.

\begin{assumption}
    \label{ass4}
    Assume that there exists constants $\gamma_{f}\in (0,1]$, $\gamma_{g}\in (0,1]$, $H_1>0$ and $H_2>0$ such that
    \begin{eqnarray*}
        &&|f(s,x)-f(t,x)|\leq H_1(1+|x|^{\a+1})(s-t)^{\gamma_{f}},\\
        &&|g(s,x)-g(t,x)|\leq H_2(1+|x|^{\a+1})(s-t)^{\gamma_{g}},
    \end{eqnarray*}
    for any $x,y \in \mathbb{R}^d$ and any $s,t\in [0,T]$.
    \par
\end{assumption}

\medskip 
\noindent
Now we introduce the procedure of constructing the Milstein-type method discussed in this paper.

\noindent
{\bf Step 1.} Based on the formates of the coefficients, we choose a strictly increasing continuous function $\m:\RR_+\rightarrow \RR_+$ such that $\m(u)\rightarrow \infty$ as $u\rightarrow \infty$ and for any $l = 1,2, ...,d$.
\begin{align*}
    \sup_{0\leq t\leq T}\sup_{|x|\leq u}(|f(t,x)|\vee|g(t,x)|\vee|G^{l}(t,x)|)\leq \m(u),\quad u\geq 1.
\end{align*}

\noindent
{\bf Step 2.} We choose a constant $\hat{\k}\geq 1\wedge \m(1)$ and a strictly decreasing function $\k:(0,1]\rightarrow [\m(1),\infty)$ such that
\begin{align}
    \label{equ0}
    h^{1/4}\k(h)\leq \hat{\k} \quad \text{for any}~ h \in (0,1] \quad \text{and} \quad  \lim_{h \rightarrow 0}\k(h)=\infty.
\end{align}

\noindent
{\bf Step 3.} Since the inverse function of $\m$, denoted by $\m^{-1}$, is a strictly increasing continuous function from $[\m(0),\infty)$ to $\RR_+$, for a given step size $h \in (0,1]$ we define the truncated mapping by
\begin{align*}
    \pi_{h}(x)=\left(|x|\wedge\m^{-1}(\k(h))\right)\frac{x}{|x|},
\end{align*}
where $x/|x|$ is set to be 0 if $x=0$. Then we define the truncated functions by
\begin{align*}
    f_{h}(t,x)=f(t,\pi_{h}(x)),\quad g_{h}(t,x)=g(t,\pi_{h}(x)),\quad G^{l}_{h}(t,x)=G^{l}(t,\pi_{h}(x)).
\end{align*}
for any $x\in \RR^d$ and $l = 1,2, ...,d$. It is not hard to see that for any $t\in[0,T]$ and any $x \in \RR^d$,
\begin{align}
    \label{equ01}
    |f_{h}(t,x)|\vee |g_{h}(t,x)|\vee |G^{l}_{h}(t,x)|\leq \m(\m^{-1}(\k(h)))=\k(h).
\end{align}

we can also obtain the fact that there exists a positive constant $\hat{M}$ such that
\begin{eqnarray*}
    &&|\frac{\partial f_{h}(t,x)}{\partial x}|\vee|\frac{\partial^{2} f_{h}(t,x)}{\partial x^{2}}|\vee|\frac{\partial g_{h}(t,x)}{\partial x}|\vee|\frac{\partial^{2} g_{h}(t,x)}{\partial x^{2}}|\leq \hat{M},
\end{eqnarray*}
for any $t\in[0,T]$ and all $x \in \RR^d$.

\noindent
{\bf Step 4.} Now we turn to discretise the process $E(t)$ in a finite time interval $[0,T]$ for any given $T>0$. For the given step size h, set $t_i = ih$ and let $\Delta_i$ be independently identically sequence satisfying $\Delta_i = D(h)$ in distribution for $i=0,1,2,...$. By the iteration, $D_h(t_i) = D_h(t_{i-1} )+ \Delta_i$ with $D_h(0) = 0$, the sample path of $D(t)$ can be simulated. And we stop the iteration for some positive integer $N$ when
\begin{align*}
    T \in [ D_h(t_{N}), D_h(t_{N+1}))
\end{align*}
holds.

\noindent
{\bf Step 5.} The discretised $E(t)$, denoted by $E_h(t)$, can be found by
\begin{align}\label{findEht}
    E_h(t) = \big(\min\{n; D_h(t_n) > t\} - 1\big)h,
\end{align}
for $t \in [0,T]$. It is not hard to see $E_h(t) = ih$ for $t \in \left[ D_h(t_{i}), D_h(t_{i+1})\right)$.

For $i = 0,1,2, ..., N$, denote $\tau _i = D_h(t_i)$. Then it can be observed that
\begin{align}\label{eq:Eh}
    E_h(\tau _i) = E_h(D_h(t_i)) = ih.
\end{align}

\noindent
{\bf Step 6.} Finally by setting $X_0=Y(0)$, the discrete version of the Milstein method is defined as
\begin{align}\label{num-ori}
    X_{\tau_{n+1}}=&X_{\tau_{n}}+f_{h}(\tau_n,X_{\tau_{n}})\bigg(E_h(\tau_{n+1}) - E_h(\tau_n)\bigg)\nonumber \\
    &+g_{h}(\tau_n,X_{\tau_{n}})\bigg(W(E_h(\tau_{n+1})) - W(E_h(\tau_n)) \bigg) \nonumber\\
    &+\dfrac{1}{2}\sum_{l=1}^{d}g^{l}_{h}(\tau_n,X_{\tau_{n}})G^{l}_{h}(\tau_n,X_{\tau_{n}})\bigg(\Delta W^{2}(E_h(\tau_{n})) -\Delta(E_h(\tau_n)) \bigg).
\end{align}
It should be noted that $\{\tau_n\}_{n=1,2,...,N}$ is a random sequence but independent from the Wiener process. In addition, it is not hard to see from \eqref{eq:Eh} that
\begin{align*}
    E_h(\tau_{n+1}) - E_h(\tau_n) = h  ~~ \text{and} ~~  W(E_h(\tau_{n+1})) - W(E_h(\tau_n)) = W((n+1)h) -W(nh).
\end{align*}

Now we present the continuous version of \eqref{num-ori}, as it is more convenient to use it in our proofs.

For any $t\in[0,T]$ and any $x \in \RR^d$, set
\begin{align*}
    Lg_h(t,x):=\sum_{l=1}^{d}g_h^{l}(t,x)G_h^{l}(t,x),
\end{align*}

For any $t\in[0, T]$. the continuous version of our Milstein method is 
\begin{align}
    \label{equ04}
        X(t)=&X(0)+\int_{0}^{t}f_{h}{(\bar{\tau}(s),\bar{X}(s))dE(s)}+\int_{0}^{t}g_{h}{(\bar{\tau}(s),\bar{X}(s))dW(E(s))}\nonumber\\
        &+\int_{0}^{t}Lg_{h}{(\bar{\tau}(s),\bar{X}(s))\Delta W(E_h(s))dW(E(s))}, 
\end{align}
where $\bar{\tau}(s)=\tau_n\ind_{[\tau_n,\tau_{n+1})}(s)$, $\bar{X}(t)=\sum_{n=0}^{N}X_{\tau_n}\ind_{[\tau_n,\tau_{n+1})}(t)$ and
\begin{align*}
\Delta W(E_h(s))=\sum_{i=1}^{N}\ind_{\left\{\tau_i\leqslant s < \tau_{i+1}\right\} }\big(W(E_h(s))-W(E_h(\tau_i))\big).
\end{align*}

\noindent
The following version of the Taylor expansion is essential for proofs in our paper.

Given $\psi  :\RR^{d+1}\to \RR^{d}$ is a third-order continuously differentiable function, for $z, z^{*}\in \RR^{d+1}$ we have 
\begin{align*}
    \begin{split}
        \psi (z)-\psi (z^{*})=\psi ^{'}(z)|_{z=z^{*}}(z-z^{*})+R_{\psi }(z,z^{*}),
    \end{split}
\end{align*}
where
\begin{align*}
    \begin{split}
        R_{\psi }(z,z^{*})=\int_{0}^{1}(1-\theta)\psi ^{''}(z)|_{z=z^{*}+\theta (z-z^{*})}(z-z^{*},z-z^{*})d\theta.
    \end{split}
\end{align*}
Here,  $\psi ^{'}$ and $\psi ^{''} $are defined in the following way, for any $z, \hat{z}, \tilde{z} \in \RR^{d+1}$.
\begin{align*}
    \psi ^{'}(z)(\hat{z})=\sum_{i=1}^{d+1}\dfrac{\partial \psi }{\partial z^{i}}\hat{z}_{i},
    \quad \psi^{''}(z)(\hat{z},h)=\sum_{i,k=1}^{d+1}\dfrac{\partial^{2} \psi  }{\partial z^{i}\partial z^{k}}\hat{z}_{i} \tilde{z}_{k},
\end{align*}
where $\psi =(\psi _{1},\psi _{2},...,\psi _{d})^{T}$, $\psi_{j}: \RR^{d+1} \to \RR$ for $j=1,2,...d$,  and $\dfrac{\partial \psi }{\partial z^{i}}=(\dfrac{\partial \psi _{1}}{\partial z^{i}},\dfrac{\partial \psi _{2}}{\partial z^{i}},...,\dfrac{\partial \psi _{d}}{\partial z^{i}})^{T}$ for $i=1,2,...,d+1$.

In the paper, we employ the Taylor expansion above by using one dimension for the time variable and $d$ dimensions for the state variable, to be more precise, we set $z=(\eta,\bar{x})$ and $z^{*}=(\eta,x^{*})$ for $\eta \in \RR_{+}$ and $\bar{x},x^{*}\in \RR^{d}$. It is clear that $z-z^{*}=(0,\bar{x}-x^{*})$, in the case. Therefore, for $\psi :\RR_{+}\times \RR^{d}\to \RR^{d}$ we have that
\begin{align*}
    \begin{split}
        \psi (\eta ,\bar{x})-\psi (\eta ,x^{*})=\psi ^{'}(\eta ,x)\big|_{x=x^{*}}(\bar{x}-x^{*})+R_{\psi }(\eta,\bar{x},x^{*}),
    \end{split}
\end{align*}
where
\begin{align*}
    \begin{split}
        R_{\psi }(\eta,\bar{x},x^{*})=\int_{0}^{1}(1-\theta )\psi ^{''}(\eta,x)|_{x=x^{*}+\theta (\bar{x}-x^{*})}(\bar{x}-x^{*},\bar{x}-x^{*})d\theta ,
    \end{split}
\end{align*}
for any $\eta \in \RR_{+}$ and $\bar{x},x^{*}\in \RR^{d}$. In this case, for any $x,\bar{j},\bar{h} \in \RR^{d},\psi ^{'}$ and $\psi ^{''}$ are defined by
\begin{align*}
    \psi ^{'}(\eta,x)(\bar{j})=\sum_{i=1}^{d}\dfrac{\partial \psi  }{\partial x^{i}}\bar{j}_{i},
    \quad\psi^{''}(\eta,x)(\bar{j},\bar{h})=\sum_{i,k=1}^{d}\dfrac{\partial^{2} \psi  }{\partial x^{i}\partial x^{k}}\bar{j}_{i}\bar{h}_{k}.
\end{align*}
respecttively. Here, $\psi =(\psi _{1},\psi _{2},...,\psi _{d})^{T}$,$\dfrac{\partial \psi }{\partial x^{i}}=(\dfrac{\partial \psi _{1}}{\partial x^{i}}, \dfrac{\partial \psi _{2}}{\partial x^{i}},...,\dfrac{\partial \psi _{d}}{\partial x^{i}})^{T}$,   $\bar{j}=(\bar{j}_{1},\bar{j}_{2},...,\bar{j}_{d})^{T}$ and $\bar{h}=(\bar{h}_{1},\bar{h}_{2},...,\bar{h}_{d})^{T}$.

Setting $\eta=\bar{\tau}(t), \bar{x}=X(t)$ and $x^{*}=\bar{X}(t)$, we derive from above, that for any fixed $t\in [0,T]$,
\begin{align}
        \label{le210}
        \psi (\bar{\tau}(t),X(t))-\psi (\bar{\tau}(t),\bar{X}(t))
        =&\psi ^{'}(\bar{\tau}(t),x)\big|_{x=\bar{X}(t)}\int_{0}^{t}g_{h}(\bar{\tau}(s),\bar{X}(s))dW(E(s))\nonumber\\
        &+ \tilde{R}_{\psi }(t,X(t),\bar{X}(t)),
\end{align}
Here
\begin{align}
        \tilde{R}_{\psi }(t,X(t),\bar{X}(t))
        =&\psi ^{'}(\bar{\tau}(t),x)\big|_{x=\bar{X}(t)}\bigg(\int_{0}^{t}f_{h}(\bar{\tau}(s),\bar{X}(s))dE(s)\nonumber\\
        &+\int_{0}^{t}Lg_{h}(\bar{\tau}(s),\bar{X}(s))\Delta W(E_{h}(s))dW(E(s))\bigg)\nonumber\\
        &+R_{\psi }(\bar{\tau}(t),X(t),\bar{X}(t)).
\end{align}
Thus, resplacing $\psi $ by $g_{h}$, we obtain
\begin{align}
        \label{the2_12}
        &\tilde{R}_{g_{h}}(t,X(t),\bar{X}(t)) \nonumber\\
        =&g_{h}(\bar{\tau}(t),X(t))-g_{h}(\bar{\tau}(t),\bar{X}(t))-Lg_{h}(\bar{\tau }(t),\bar{X}(t))\Delta W(E_{h}(t)).
\end{align}

At the end of this section, we mention some known results. For the proofs of Lemmas \ref{lemma2-6} and \ref{lemma23}, we refer the readers to \cite{HLM2018}. The proof of Lemma \ref{lemma_E} can be found in \cite{JK2016}. Lemma \ref{lemY} is borrowed from \cite{LLTarXiv}.

\begin{lemma}
    \label{lemma2-6}
    Let Assumption \ref{ass1} hold. For all $h\in (0,1]$
    \begin{align*}
        &|f_{h}(t,x)-f_{h}(t,y)|\vee |g_{h}(t,x)-g_{h}(t,y)|\vee |Lg_{h}(t,x)-Lg_{h}(t,y)|\\
        \leq& C(1+|x|^{\a}+|y|^{\a})|x-y|
    \end{align*}
holds for all $t\in(0,T]$ and $x,y \in \RR^{d}$.
\end{lemma}

\begin{lemma}
    \label{lemma23}
    Let Assumption \ref{ass3} hold. Then, for all $h\in (0,1]$, we have
    \begin{align*}
        x^{\mathrm{T}}f_{h}(t,x) +(5q-1)|g_{h}(t,x)|^2 \leq \hat{K}_1(1+|x|^2),\quad \forall x\in \RR^d,
    \end{align*}
    where $\hat{K}_1=2K_1\left(1\vee \frac{1}{\m^{-1}(\k(1))}\right)$.
\end{lemma}

\begin{lemma}
    \label{lemma_E}
    For any $t_i\leq t\leq t_{i+1}$, there exists a constant $c$ such that
    \begin{align*}
        |E_h(t)-E_h(t_i)|\leq |E_h(t_{i+1})-E_h(t_i)|\leq ch.
    \end{align*}
\end{lemma}

\begin{lemma}
    \label{lemY}
    Suppose Assumption \ref{ass1} and \ref{ass3} hold. Then, for any $p\in [2,q)$
    \begin{align*}
        \EE \left( \sup _{0\leq t\leq T}|Y(t)|^p \right)<\infty.
    \end{align*}
\end{lemma}

Briefly speaking, Lemmas \ref{lemma2-6} and \ref{lemma23} indicate that, to some extended, the truncated functions $f_h$ and $g_h$ inherit Assumptions  \ref{ass1} and \ref{ass3}. Lemma \ref{lemma_E} is useful for the analysis of the convergence order of $E_h(t)$. Lemma \ref{lemY} states the moment boundedness of the underlying solution.

\section{Lemmas prepared for main results}

Lemmas that will be used in the proofs of main results in Section 4 are presented and proved in this section. 

\begin{lemma}
    \label{lemma2}
    For any $h\in (0,1]$ and any $\hat{p}>2$, we have
    \begin{align}
        \label{equ06}
        \EE_W|X(t)-\bar{X}(t)|^{\hat{p}}\leq c_{\hat{p}}h^{\hat{p}/2}\left(\k(h)\right)^{\hat{p}},\quad \forall t\geq 0,
    \end{align}
    where $c_{\hat{p}}=c\left(\frac{\hat{p}(\hat{p}-1)}{2}\right)^{\frac{\hat{p}}{2}}3^{\hat{p}-1}$, consequently,
    \begin{align}
        \label{equ07}
        \lim_{h\rightarrow 0}\EE_W|X(t)-\bar{X}(t)|^{\hat{p}}=0,\quad \forall t\geq 0.
    \end{align}
\end{lemma}
\noindent
{\bf Proof.}
    Fix any $h\in (0,1]$, $\hat{p}>2$ and $t\geq 0$. There is a unique integer $n\geq 0$ such that $\tau_n\leq t < \tau_{n+1}$. By properties of the basic inequality, we then derive from \eqref{equ04} that
    \begin{align}
        \label{lem250}
            &|X(t)-\bar{X}(t)|^{\hat{p}}\nonumber\\
            =&|X(t)-X(\tau_n)|^{\hat{p}}\nonumber\\
            =&\bigg|\int_{\tau_n}^{t}{f_{h}(\bar{\tau}(s),\bar{X}(s))dE(s)} +\int_{\tau_n}^{t}{g_{h}(\bar{\tau}(s),\bar{X}(s))dW(E(s))} \nonumber\\
            &+\int_{\tau_n}^{t}{Lg_{h}(\bar{\tau}(s),\bar{X}(s))\Delta W(E_{h}(s))dW(E(s))}\bigg|^{\hat{p}}\nonumber\\
            \leq&3^{\hat{p}-1}\bigg(\left|\int_{\tau_n}^{t}{f_{h}(\bar{\tau}(s),\bar{X}(s))dE(s)} \right|^{\hat{p}}+\left|\int_{\tau_n}^{t}{g_{h}(\bar{\tau}(s),\bar{X}(s))dW(E(s))} \right|^{\hat{p}}\nonumber\\
            &+\left|\int_{\tau_n}^{t}{Lg_{h}(\bar{\tau}(s),\bar{X}(s))\Delta W(E_{h}(s))dW(E(s))} \right|^{\hat{p}}\bigg).
    \end{align}
Now, we estimate those three terms inside the bracket on the right hand side of the last inequality of \eqref{lem250}
By the H\"older inequality, the first term can be estimated by
    \begin{align}
        \label{lem38}
        \EE_W\left|\int_{\tau_n}^{t}{f_{h}(\bar{\tau}(s),\bar{X}(s))dE(s)} \right|^{\hat{p}}\leq h^{\hat{p}-1}\EE_W\int_{\tau_n}^{t}{\left|f_{h}(\bar{\tau}(s),\bar{X}(s))\right|^{\hat{p}} dE(s)}.
    \end{align}

The second item 
 Let  $x(t)=\int_{\tau_n}^{t}{g_{h}(\bar{\tau}(s),\bar{X}(s))dW(E(s))}$, so we have
    \begin{align*}
            &\EE_W\left| x(t)\right|^{\hat{p}}\\
            =&\dfrac{\hat{p}}{2}\EE_W\int_{\tau_n}^{t}\left(\left| x(s)\right|^{\hat{p}-2}\big|{g_{h}(\bar{\tau}(s),\bar{X}(s))\big|^{2}+({\hat{p}-2})\left| x(s)\right|^{\hat{p}-4}}\big|x^{T}(s)g(s)\big|^{2}\right)dE(s)\nonumber\\ 
            \leq& \frac{\hat{p}(\hat{p}-1)}{2}\EE_W\int_{\tau_n}^{t}\left| x(s)\right|^{\hat{p}-2}\left|g_{h}(\bar{\tau}(s),\bar{X}(s))\right|^2dE(s)\\
            \leq& \frac{\hat{p}(\hat{p}-1)}{2}\left(\EE_W\int_{\tau_n}^{t}\left| x(s)\right|^{\hat{p}}dE(s)\right)^{\frac{\hat{p}-2}{\hat{p}}}\left(\EE_W\int_{\tau_n}^{t}\left|g_{h}(\bar{\r}(s),\bar{X}(s))\right|^{\hat{p}}dE(s)\right)^{\frac{2}{\hat{p}}}\\
            =&  \frac{\hat{p}(\hat{p}-1)}{2}\left(\int_{\tau_n}^{t}\EE_W\left| x(s)\right|^{\hat{p}}dE(s)\right)^{\frac{\hat{p}-2}{\hat{p}}}\left(\EE_W\int_{\tau_n}^{t}\left|g_{h}(\bar{\r}(s),\bar{X}(s))\right|^{\hat{p}}dE(s)\right)^{\frac{2}{\hat{p}}}.
\end{align*}

    Note that $\EE_W\left| x(t)\right|^{\hat{p}}$ is nondecreasing in $t$, it then follows\\
    \begin{align*}
        \begin{split}
            \EE_W\left| x(t)\right|^{\hat{p}}
            &\leq \frac{\hat{p}(\hat{p}-1)}{2}\left[ch\EE_W\left| x(t)\right|^{\hat{p}}\right]^{\frac{\hat{p}-2}{\hat{p}}}\left(\EE_W\int_{\tau_n}^{t}\left|g_{h}(\bar{\tau}(s),\bar{X}(s))\right|^{\hat{p}}dE(s)\right)^{\frac{2}{\hat{p}}}.
        \end{split}
    \end{align*}
It is obtained by further shifting and simplification,
    \begin{align*}
       \begin{split}
            \EE_W\left| x(t)\right|^{\hat{p}}
            \leq \left(\frac{\hat{p}(\hat{p}-1)}{2}\right)^{\frac{\hat{p}}{2}}ch^{\frac{\hat{p}-2}{2}}\EE_W\int_{\tau_n}^{t}\left|g_{h}(\bar{\tau}(s),\bar{X}(s))\right|^{\hat{p}}dE(s).\\
        \end{split}
    \end{align*}
    So we can have
    \begin{align}
            \label{lem39}
            &\EE_W\left|\int_{\tau_n}^{t}{g_{h}(\bar{\tau}(s),\bar{X}(s)) dW(E(s))}\right|^{\hat{p}}\nonumber\\
            \leq& \left(\frac{\hat{p}(\hat{p}-1)}{2}\right)^{\frac{\hat{p}}{2}}ch^{\frac{\hat{p}-2}{2}}\EE_W\int_{\tau_n}^{t}\left|g_{h}(\bar{\tau}(s),\bar{X}(s))\right|^{\hat{p}}dE(s).
    \end{align}

    The third item in the above brackets, using the same way as the second item, we can see
    \begin{align}
            \label{lem310}
            &\EE_W\left|\int_{\tau_n}^{t}{Lg_{h}(\bar{\tau}(s),\bar{X}(s))\Delta W(E(s))dW(E_{h}(s))}\right|^{\hat{p}}\nonumber\\
            &\leq \left(\frac{\hat{p}(\hat{p}-1)}{2}\right)^{\frac{\hat{p}}{2}}ch^{\frac{\hat{p}-2}{2}}\EE_W\int_{\tau_n}^{t}\left|Lg_{h}(\bar{\tau}(s),\bar{X}(s))\Delta W(E_{h}(s))\right|^{\hat{p}}dE(s).
    \end{align}
    Substituting the estimates \eqref{lem38} ,\eqref{lem39} and \eqref{lem310}  into \eqref{lem250}, we use the \eqref{eq:Eh} and lemma \ref{lemma_E} and \eqref{equ01}, we have
    \begin{align*}
            \EE_W|X(t)-\bar{X}(t)|^{\hat{p}} &\leq
            3^{\hat{p}-1}\bigg(h^{\hat{p}-1}\EE_W\int_{\tau_n}^{t}{\left|f_{h}(\bar{\tau}(s),\bar{X}(s))\right|^{\hat{p}}dE(s)}\\
            &+\left(\frac{\hat{p}(\hat{p}-1)}{2}\right)^{\frac{\hat{p}}{2}}ch^{\frac{\hat{p}-2}{2}}\EE_W\int_{\tau_n}^{t}\left|g_{h}(\bar{\tau}(s),\bar{X}(s))\right|^{\hat{p}}dE(s)\\
            &+\left(\frac{\hat{p}(\hat{p}-1)}{2}\right)^{\frac{\hat{p}}{2}}ch^{\frac{\hat{p}-2}{2}}\EE_W\int_{\tau_n}^{t}\bigg|Lg_{h}(\bar{\tau}(s),\bar{X}(s))\\
            &\times \Delta W(E_{h}(s))\bigg|^{\hat{p}}dE(s)\bigg)\\
            &\leq3^{\hat{p}-1}\bigg(h^{\hat{p}-1}ch(\k(h))^{\hat{p}}+\left(\frac{\hat{p}(\hat{p}-1)}{2}\right)^{\frac{\hat{p}}{2}}ch^{\frac{\hat{p}-2}{2}}h(\k(h))^{\hat{p}}\\
            &+\left(\frac{\hat{p}(\hat{p}-1)}{2}\right)^{\frac{\hat{p}}{2}}ch^{\frac{\hat{p}-2}{2}}h^{\frac{\hat{p}}{2}}h(\k(h))^{2\hat{p}}\bigg)\\
            &\leq c_{\hat{p}}\left(h^{\hat{p}-1}h(\k(h))^{\hat{p}}+h^{\hat{p}/2-1}h(\k(h))^{\hat{p}}+h^{\hat{p}/2}h^{\hat{p}/2}(\k(h))^{2\hat{p}}\right)\\
            &\leq c_{\hat{p}}\left(h^{\hat{p}}(\k(h))^{\hat{p}}+h^{\hat{p}/2}(\k(h))^{\hat{p}}+h^{\hat{p}}(\k(h))^{2\hat{p}}\right)\\
            &\leq c_{\hat{p}}h^{\hat{p}/2}(\k(h))^{\hat{p}},
    \end{align*}
    where $c_{\hat{p}}=c\left(\frac{\hat{p}(\hat{p}-1)}{2}\right)^{\frac{\hat{p}}{2}}3^{\hat{p}-1}$, this completes the proof of \eqref{equ06}. Noting from \eqref{equ0}, we have $h^{\hat{p}/2}(\k(h))^{\hat{p}}\leq h^{\hat{p}/4}$. Then, \eqref{equ07} can be derived from \eqref{equ06}.
\eproof

Now, we prove the boundedness of the pth  moment the numerical solution.
\begin{lemma}
    \label{lemma3}
    Let Assumptions \ref{ass1} and \ref{ass3} hold. Then
    \begin{align}
        \label{equ08}
        \sup_{0< h \leq 1}\EE [\sup_{0\leq t\leq T}|X(t)|^p]\leq C, \quad \forall T>0,
    \end{align}
    where $C=\bigg(2|X(0)|^p+4c_{p}^{\frac{1}{2}}\hat{k}E(t)+2(5p^{2}-p)c\hat{k}E(t)\bigg)e^{3(2p\hat{K}_1\vee 2(p-2)\vee2(5p^{2}-p))E(T)}$ is a constant dependent on $X(0)$, $p$, $T$,$c_{p}$,$\hat{k}$ and  $\hat{K}_1$, but independent from $h$.
\end{lemma}
{\bf Proof.}
    Define the stopping time $\zeta_{\ell}:=\inf\{ t\geq 0; |X(t)|>\ell \}$ for some positive integer $\ell$. It can be seen that
    \begin{align*}
        \int_{0}^{t}{\EE_W\left(\sup_{0\leq s \leq t\wedge \zeta_{\ell}}|X(s)|^p\right)dE(r)}\leq \ell^pE(t).
    \end{align*}
    Fix any $h\in(0,1]$ and $T\geq 0$. By the It\^o formula, we derive from \eqref{equ04} that, for $0\leq u\leq t\wedge \zeta_{\ell}$,
    \begin{align}
        |X(u)|^p=|X(0)|^p+A_u+M_u,
    \end{align}
    where
    \begin{align*}
            A_u:=&\int_{0}^{u}\bigg( p|X(s)|^{p-2}X^{\mathrm{T}}(s)f_{h}(\bar{\tau}(s),\bar{X}(s))+\frac{1}{2}p(p-1)|X(s)|^{p-2}|g_{h}(\bar{\tau}(s),\bar{X}(s))\\
            &+Lg_{h}(\bar{\tau}(s),\bar{X}(s))\Delta W(E_{h}(s))|^2\bigg)dE(s),
    \end{align*}
    \begin{align*}
            M_u:=
            \int_{0}^{u}{p|X(s)|^{p-1}|g_{h}(\bar{\tau}(s),\bar{X}(s))+Lg_{h}(\bar{\tau}(s),\bar{X}(s))\Delta W(E_{h}(s))|dW(E(s))}.
    \end{align*}
    It can be noted that the stochastic integral $(M_u)_{u\geq 0}$ is a local martingale with quadratic variation
    \begin{align*}
        [M,M]_u=\int_{0}^{u}{p^2|X(s)|^{2p-2}\big|g_{h}(\bar{\tau}(s),\bar{X}(s))+Lg_{h}(\bar{\tau}(s),\bar{X}(s))\Delta W(E_{h}(s))\big|^2 dE(s)},
    \end{align*}
    For $0\leq s\leq t\wedge \zeta_{\ell}$,
    \begin{align*}
            &p^2|X(s)|^{2p-2}\big|g_{h}(\bar{\tau}(s),\bar{X}(s))+Lg_{h}(\bar{\tau}(s),\bar{X}(s))\Delta W(E_{h}(s))\big|^2\\
            &\leq p^2|X(s)|^p|X(s)|^{p-2}\big|g_{h}(\bar{\tau}(s),\bar{X}(s))+Lg_{h}(\bar{\tau}(s),\bar{X}(s))\Delta W(E_{h}(s))\big|^2\\
            &\leq p^2\left(\sup_{0\leq u\leq t\wedge \zeta_{\ell}}|X(u)|^p\right)\big|X(s)|^{p-2}|g_{h}(\bar{\tau}(s),\bar{X}(s))+Lg_{h}(\bar{\tau}(s),\bar{X}(s))\Delta W(E_{h}(s))\big|^2.
    \end{align*}
    By using the inequality $(ab)^{1/2}\leq a/\l +\l b$ valid for any $a,b\geq 0$ and $\l > 0$, we can see that for $0\leq u \leq t\wedge \zeta_{\ell}$,
    \begin{align*}
            &([M,M]_u)^{1/2}\\
            \leq& p\bigg(\sup_{0\leq u\leq t\wedge \zeta_{\ell}}|X(u)|^p\int_{0}^{u}|X(s)|^{p-2}\big|g_{h}(\bar{\tau}(s),\bar{X}(s))\\
            &+Lg_{h}(\bar{\tau}(s),\bar{X}(s))\Delta W(E_{h}(s))\big|^2 dE(s)\bigg)^{1/2}\\
            \leq& p\bigg(\frac{\sup_{0\leq u\leq t\wedge \zeta_{\ell}}|X(u)|^p}{2p}+2p\int_{0}^{u}|X(s)|^{p-2}\big|g_{h}(\bar{\tau}(s),\bar{X}(s))\\
            &+Lg_{h}(\bar{\tau}(s),\bar{X}(s))\Delta W(E_{h}(s))\big|^2 dE(s)\bigg).
    \end{align*}
    We have expectations for $A_u$ and $M_u$, respectively
    \begin{align}
            \EE_W(A_u)=&\EE_W\bigg(\sup_{0\leq u\leq t\wedge \zeta_{\ell}}\int_{0}^{u}\bigg(p|X(s)|^{p-2}X^{\mathrm{T}}(s)f_{h}(\bar{\tau}(s),\bar{X}(s))\nonumber\\ &+\frac{1}{2}p(p-1)|X(s)|^{p-2}\big|g_{h}(\bar{\tau}(s),\bar{X}(s))\nonumber\\
            &+Lg_{h}(\bar{\tau}(s),\bar{X}(s))\Delta W(E_{h}(s))\big|^2\bigg) dE(s)\bigg),
    \end{align}
\begin{align}
\EE_W(M_u)=&\EE_W\bigg(\frac{1}{2}\sup_{0\leq u\leq t\wedge\zeta_{\ell}}|X(u)|^p+\sup_{0\leq u\leq t\wedge\zeta_{\ell}}\int_{0}^{u}2p^2|X(s)|^{p-2}\big|g_{h}(\bar{\tau}(s),\bar{X}(s))  \nonumber\\ 
 &+Lg_{h}(\bar{\tau}(s),\bar{X}(s))\Delta W(E_{h}(s))\big|^2 dE(s)\bigg).
\end{align}
    Take the expectation for (3.12), then substitute (3.13) and (3.14), and use the basic inequality $(a+b)^{2}\leq2(a^{2}+b^{2})$, we can have
\begin{align}
&\EE_W\bigg(\sup_{0\leq u\leq t\wedge \zeta_{\ell}}|X(u)|^p\bigg)  \nonumber\\ 
=&|X(0)|^p+\EE_W(A_u)+\EE_W(M_u) \nonumber\\ \nonumber
\leq& |X(0)|^p+\frac{1}{2}\EE_W\bigg(\sup_{0\leq u\leq t\wedge \zeta_{\ell}}|X(u)|^p\bigg) \\ \nonumber
\quad&+\EE_W\bigg(\sup_{0\leq u\leq t\wedge \zeta_{\ell}}\int_{0}^{u}p |X(s)|^{p-2}\bigg(X^{\mathrm{T}}(s)f_{h}(\bar{\tau}(s),\bar{X}(s)) \\ \nonumber
 \quad&+(5p-1)|g_{h}(\bar{\tau}(s),\bar{X}(s))|^2\bigg)dE(s)\bigg)
+\big(p(p-1)+4p^{2}\big) \\ \nonumber
\quad& \times \EE_W\bigg(\sup_{0\leq u\leq t\wedge\zeta_{\ell}}\int_{0}^{u}|X(s)|^{p-2}|Lg_{h}(\bar{\tau}(s),\bar{X}(s)) \Delta W(E_{h}(s))|^2dE(s)\bigg) \\ \nonumber
\leq& |X(0)|^p+\frac{1}{2}\EE_W\bigg(\sup_{0\leq u\leq t\wedge \zeta_{\ell}}|X(u)|^p\bigg) \\ \nonumber
\quad&+\EE_W\bigg(\sup_{0\leq u\leq t\wedge \zeta_{\ell}}\int_{0}^{u}p|X(s)|^{p-2} \bigg(\bar{X}^{\mathrm{T}}(s)f_{h}(\bar{\tau}(s),\bar{X}(s)) \\ \nonumber
\quad&+(5p-1)|g_{h}(\bar{\tau}(s),\bar{X}(s))|^2\bigg)dE(s)\bigg) \\ \nonumber
\quad& +\EE_W\bigg(\sup_{0\leq u\leq t\wedge \zeta_{\ell}}\int_{0}^{u}p|X(s)|^{p-2}(X(s)-\bar{X}(s))^{\mathrm{T}}f_{h}(\bar{\tau}(s),\bar{X}(s))dE(s)\bigg) \\ \nonumber
\quad& 
+(5p^{2}-p) \EE_W\bigg(\sup_{0\leq u\leq t\wedge \zeta_{\ell}}\int_{0}^{u}|X(s)|^{p-2}|Lg_{h}(\bar{\tau}(s),\bar{X}(s))\Delta W(E_{h}(s))|^2dE(s)\bigg).\nonumber \\ 
\end{align}
    Therefore, for any $0\leq u\leq t\wedge \zeta_{\ell}$, by Lemma \ref{lemma23} and the Young inequality
    \begin{align*}
        a^{p-2}b\leq \frac{p-2}{p}a^p + \frac{2}{p}b^{p/2},\quad \forall a,b \geq 0.
    \end{align*}
    we can get from (3.15)
    \begin{align*}
        \label{e316}
            &\EE_W\bigg(\sup_{0\leq u\leq t\wedge \zeta_{\ell}}|X(u)|^p\bigg)\\
            \leq& |X(0)|^p+\frac{1}{2}\EE_W\bigg(\sup_{0\leq u\leq t\wedge \zeta_{\ell}}|X(u)|^p\bigg)\\
            \quad&+p\hat{K}_1\EE_W\left(\sup_{0\leq u\leq t\wedge \zeta_{\ell}}\int_{0}^{u}{|X(s)|^{p-2}(1+|\bar{X}(s)|^2)dE(s)}\right)\\
            \quad& +\EE_W\sup_{0\leq u\leq t\wedge \zeta_{\ell}}\bigg((p-2)\int_{0}^{u}{|X(s)|^p dE(s)}\\
            \quad&+2\int_{0}^{u}{|X(s)-\bar{X}(s)|^{p/2}|f_{h}(\bar{\tau}(s),\bar{X}(s))|^{p/2}dE(s)}\bigg)\\
            \quad& +(5p^{2}-p)\EE_W\bigg(\sup_{0\leq u\leq t\wedge \zeta_{\ell}}\int_{0}^{u}{|X(s)|^{p-2}|Lg_{h}(\bar{\tau}(s),\bar{X}(s))\Delta W(E_{h}(s))|^2dE(s)}\bigg).\\
    \end{align*}
    Thus, for any $0\leq u\leq t\wedge \zeta_{\ell}$ and apply basic inequality, we have
    \begin{align}
            \EE_W\left(\sup_{0\leq u\leq t\wedge \zeta_{\ell}}|X(t)|^p\right)&\leq 2|X(0)|^p+2p\hat{K}_1\EE_W\int_{0}^{t}|X(t\wedge \zeta_{\ell})|^{p-2} \nonumber\\
            \quad& \times(1+|\bar{X}(s)|^2)dE(s) \nonumber\\
            \quad& +2(p-2)\int_{0}^{t}\EE_W |X(t\wedge \zeta_{\ell})|^p dE(s)+I_{1}+I_{2},
    \end{align}
    where
    \begin{align*}
            I_{1}=4\EE_W\int_{0}^{t}{|X(s)-\bar{X}(s)|^{p/2}|f_{h}(\bar{\tau}(s),\bar{X}(s))|^{p/2}dE(s)},
    \end{align*}
    \begin{align*}
            I_{2}=2(5p^{2}-p)\EE_W\left(\int_{0}^{t}{|X(t\wedge \zeta_{\ell})|^{p-2}|Lg_{h}(\bar{\tau}(s),\bar{X}(s))\Delta W(E_{h}(s))|^2dE(s)}\right).
    \end{align*}
    Now we deal with the $I_{1}$ item above,  by Lemma \ref{lemma2}, inequalities \eqref{equ0} and \eqref{equ01}, we have
    \begin{align*}
        \label{}
            I_{1}=&4\EE_W\int_{0}^{t}{|X(s)-\bar{X}(s)|^{p/2}|f_{h}(\bar{\tau}(s),\bar{X}(s))|^{p/2}dE(s)}\\
            \leq& 4\left(\k(h)\right)^{p/2}\int_{0}^{t}{\EE_W|X(s)-\bar{X}(s)|^{p/2}dE(s)}\\
            \leq& 4\left(\k(h)\right)^{p/2}\int_{0}^{t}{\big(\EE_W|X(s)-\bar{X}(s)|^p\big)^{1/2}dE(s)}\\
            \leq& 4c_{p}^{\frac{1}{2}}\left(\k(h)\right)^ph^{p/4}E(t)\\
            \leq& 4c_{p}^{\frac{1}{2}}\hat{k}E(t).
    \end{align*}
    We deal with the $I_{2}$ item above, by inequalities \eqref{equ01},\eqref{equ0} and lemma \ref{lemma_E}, we have
    \begin{align*}
            I_{2}&=2(5p^{2}-p)\EE_W\left(\int_{0}^{t}{|X(t\wedge \zeta_{\ell})|^{p-2}|Lg_{h}(\bar{\tau}(s),\bar{X}(s))\Delta W(E_{h}(s))|^2dE(s)}\right)\\
            &\leq2(5p^{2}-p)\EE_W\left(\int_{0}^{t}{|X(t\wedge \zeta_{\ell})|^{p-2}ch|\k(h)|^4 dE(s)}\right)\\
            &\leq2(5p^{2}-p)\left(\EE_W\int_{0}^{t}{\frac{p-2}{p}|X(t\wedge \zeta_{\ell})|^{p}dE(s)}+\EE_W\int_{0}^{u}{\frac{2}{p}|\k(h)|^{2p}ch^{\frac{p}{2}} dE(s)}\right)\\
            &\leq2(5p^{2}-p)\left(\EE_W\int_{0}^{t}{|X(t\wedge \zeta_{\ell})|^{p}dE(s)}\right)+2(5p^{2}-p)|\k(h)|^{2p}ch^{\frac{p}{2}}E(t)\\
            &\leq2(5p^{2}-p)\left(\EE_W\int_{0}^{t}{|X(t\wedge \zeta_{\ell})|^{p}dE(s)}\right)+2(5p^{2}-p)c\hat{k}E(t).\\
    \end{align*}
    We can obtain by substituing $I_{1}$ and $I_{2}$ into the \eqref{e316} \\
    \begin{align*}
            \EE_W\left(\sup_{0\leq u\leq t\wedge \zeta_{\ell}}|X(u)|^p\right)&\leq
            2|X(0)|^p+2p\hat{K}_1\EE_W\int_{0}^{t}{|X(t\wedge \zeta_{\ell})|^{p-2}(1+|\bar{X}(s)|^2)dE(s)}\\
            &\quad +2(p-2)\int_{0}^{t}{\EE_W |X(t\wedge \zeta_{\ell})|^p dE(s)}
            +4c_{p}^{\frac{1}{2}}\hat{k}E(t)\\
            &\quad+2(5p^{2}-p)\left(\EE_W\int_{0}^{t}{|X(t\wedge \zeta_{\ell})|^{p}dE(s)}\right)+2(5p^{2}-1)c\hat{k}E(t)\\
            &\leq C_1+3C_2\int_{0}^{t}{\EE_W\left(\sup_{0\leq u\leq t\wedge \zeta_{\ell}}|X(u)|^p\right) dE(s)},
    \end{align*}
    where $C_1=2|X(0)|^p+4c_{p}^{\frac{1}{2}}\hat{k}E(t)+2(5p^{2}-p)c\hat{k}E(t)$ and $C_2=2p\hat{K}_1\vee 2(p-2)\vee2(5p^{2}-p)$,
    applying the well-known Gronwall-type inequality, for any $t\in[0,T]$,
    \begin{align*}
        \EE_W\left(\sup_{0\leq u\leq t\wedge \zeta_{\ell}}|X(u)|^p\right)\leq C_1e^{(3C_2)E(t)}.
    \end{align*}
    Since $\zeta_{\ell}\rightarrow \infty$ as $\ell\rightarrow \infty$. Setting $t=T$ and letting $\ell\rightarrow \infty$ give
    \begin{align*}
        \EE_W\left(\sup_{0\leq t\leq T}|X(t)|^p\right)\leq C_1e^{(3C_2)E(T)}.
    \end{align*}
    Taking $\EE_D$ on both sides, and using the fact that $\EE_D\left(E(T)e^{E(T)}\right)<\EE_D\left(e^{2 E(T)}\right)<\EE_D\left(e^{3 E(T)}\right)< \infty$ yield,
    \begin{align*}
        \EE\left(\sup_{0\leq t\leq T}|X(t)|^p\right) \leq C,
    \end{align*}
    where $C=\bigg(2|X(0)|^p+4c_{p}^{\frac{1}{2}}\hat{k}E(t)+2(5p^{2}-p)c\hat{k}E(t)\bigg)e^{3(2p\hat{K}_1\vee 2(p-2)\vee2(5p^{2}-p))E(T)}$, as this holds for any $h\in(0,1]$ and $C$ is independent of $h$, we see the required assertion \eqref{equ08}.  
\eproof

\begin{lemma}
    \label{lemma4}
    Let Assumptions \ref{ass1},\ref{ass3},\ref{ass4} and \ref{ass5} hold, and assume that $q\geq 2(\a+1)p$ for a constant $p>2$,  then for any $\bar{p}\in [2,p)$ and $h \in (0,1]$,
    \begin{align*}
        \sup_{0<h\leq 1}\sup_{0\leq t\leq T}\left [\EE| f^{'}(t,x)|_{x=X(t)}|^{\bar{p}}\vee\EE| g^{'}(t,x)|_{x=X(t)}|^{\bar{p}}  \right ]< \infty ,
    \end{align*}
    where $f^{'}$ and $g^{'}$ denote the first partial derivatives of $f$ and $g$ with respect to the state variable $x$,respectively.
    
    We can derive it from Assumption \ref{ass5} and lemma \ref{lemma4}.
\end{lemma}  

\begin{lemma}
    \label{lemma5}
    Let Assumptions \ref{ass1},\ref{ass2},\ref{ass3}, \ref{ass4}  and \ref{ass5} hold, and assume that $q\geq 2(\a+1)p$ for a constant $p>2$, then for any $\bar{p}\in [2,p)$ and $h \in (0,1]$, $t\in [0,T]$,
    \begin{align*}
        \EE|\tilde{R}_{f}(t,X(t),\bar{X}(t))|^{\bar{p}}\vee\EE|\tilde{R}_{g}(t,X(t),\bar{X}(t))|^{\bar{p}}\vee\EE|\tilde{R}_{g_h}(t,X(t),\bar{X}(t))|^{\bar{p}} < Ch^{\bar{p}}(\k(h))^{2\bar{p}},
    \end{align*}
    where $C$ is a positive constant independent of $h$ and $t$.
\end{lemma}  
{\bf Proof.} First, for all  $0\leq t\leq T$, we give an estimate on$|R_{f}(t,X(t),\bar{X}(t))|^{\bar{p}}$ by Assusmption \ref{ass5}, lemma \ref{lemma2} and lemma \ref{lemma3}, there exists a constant C such that, we apply H{\"o}lder inequality and Jesen's inequality.
\begin{align}
        &\EE_{W}|R_{f}(t,X(t),\bar{X}(t))|^{\bar{p}}\nonumber\\
        \leq&\int^{1}_{0}(1-\theta)^{\bar{p}}\EE_{W}\big|f^{''}(\bar{\tau}(t),x)|_{x=\bar{X}(t)+\theta(X(t)-\bar{X}(t))}\nonumber\\
        &\times\big(X(t)-\bar{X}(t),X(t)-\bar{X}(t)\big)\big|^{\bar{p}}d\theta \nonumber\\
        \leq&\int^{1}_{0}\big[\EE_{W}\big|f^{''}(\bar{\tau}(t),x)|_{x=\bar{X}(t)+\theta(X(t)-\bar{X}(t))}\big|^{2\bar{p}}\EE_{W}|X(t)-\bar{X}(t)|^{4\bar{p}}\big]^{\frac{1}{2}}d\theta \nonumber\\
        \leq&C\big(1+\EE_{W}|X(t)|^{2(1+\a)\bar{p}}+\EE_{W}|\bar{X}(t)|^{2(1+\a)\bar{p}}\big)^{\frac{1}{2}}\big(\EE_{W}|X(t)-\bar{X}(t)|^{4\bar{p}}\big)^{\frac{1}{2}} \nonumber\\
        \leq&Ch^{\bar{p}}k(h)^{2\bar{p}}.
\end{align}
Then we can observe from (2.11), and the H{\"o}lder inequality that.
\begin{align}
        &\EE_{W}|\tilde{R}_{f}(t,X(t),\bar{X}(t))|^{\bar{p}}\nonumber\\
        \leq&
        C\big[h^{\bar{p}}\EE_{W}\big|f^{'}(\bar{\tau}(t),x)|_{x=\bar{X}(t)}f_{h}(\bar{\tau}(t),\bar{X}(t))\big|^{\bar{p}}\nonumber\\
        &+\dfrac{1}{2}\EE_{W}\big|f^{'}(\bar{\tau}(t),x)|_{x=\bar{X}(t)}Lg_{h}(\bar{\tau}(t),\bar{X}(t))(\Delta W(E_{h}(t))^{2}-h)\big|^{\bar{p}}\nonumber\\
        &+\EE_{W}|R_{f}(t,X(t),\bar{X}(t))|^{\bar{p}}\big]\nonumber\\
        \leq&
        C\big[h^{\bar{p}}\EE_{W}|f^{'}(\bar{\tau}(t),x)|_{x=\bar{X}(t)}f_{h}(\bar{\tau}(t),\bar{X}(t))|^{\bar{p}}\nonumber\\
        &+\dfrac{1}{2}\big(\EE_{W}|f^{'}(\bar{\tau}(t),x)|_{x=\bar{X}(t)}Lg_{h}(\bar{\tau}(t),\bar{X}(t))|^{2\bar{p}}\nonumber\EE_{W}|\Delta W(E(t))^{2}-h|^{2\bar{p}}\big)^{\frac{1}{2}}\\
        &+\EE_{W}|R_{f}(t,X(t),\bar{X}(t))|^{\bar{p}}\big].
\end{align}
We can derive from the elementary inequality $|\sum_{i=1}^{m}a_{i}|\leq m^{p-1}\sum_{i=1}^{m}|a_{i}|^{p}$ and Lemma \ref{lemma_E} that
\begin{align}
        \EE_{W}|\Delta W(E(t))^{2}-h|^{2\bar{p}}
        \leq&
        2^{2\bar{p}-1}(\EE_{W}|\Delta W(E(t))|^{4\bar{p}}+h^{2\bar{p}})\nonumber\\
        \leq&
        2^{2\bar{p}-1}(\Delta (E(t))^{2\bar{p}}+h^{2\bar{p}})\nonumber\\
        \leq&
        2^{2\bar{p}-1}(2h^{2\bar{p}})\nonumber\\
        \leq&
        2^{2\bar{p}}ch^{2\bar{p}}.
\end{align}
By using \eqref{equ01} and lemma \ref{lemma4}, we can see that for $0\leq t\leq T$,
\begin{align}
        \EE_{W}\big|f^{'}(\bar{\tau}(t),x)|_{x=\bar{X}(t)}f_{h}(\bar{\tau}(t),\bar{X}(t))\big|^{\bar{p}}\leq C(k(h))^{\bar{p}},
\end{align}
\begin{align}
        \EE_{W}\big|f^{'}(\bar{\tau}(t),x)|_{x=\bar{X}(t)}Lg_{h}(\bar{\tau}(t),\bar{X}(t))\big|^{2\bar{p}}\leq C(k(h))^{4\bar{p}}.
\end{align}
Subsittuting (3.17),(3.19),(3.20) and (3.21) into (3.18) and using the independence between $ \bar{X}(t)$and $\Delta W(t)$, we have
\begin{align*}
        \EE_{W}|\tilde{R}_{f}(t,X(t),\bar{X}(t))|^{\bar{p}}\leq Ch^{\bar{p}}(k(h))^{2\bar{p}}.
\end{align*}
Taking $\EE_D$ on the both sides, we have\\
\begin{align*}
        \EE|\tilde{R}_{f}(t,X(t),\bar{X}(t))|^{\bar{p}}\leq Ch^{\bar{p}}(k(h))^{2\bar{p}}.
\end{align*}
Similarly, we can show
\begin{align*}
    \begin{split}
        \EE|\tilde{R}_{g}(t,X(t),\bar{X}(t))|^{\bar{p}}\vee\EE|\tilde{R}_{g_{h}}(t,X(t),\bar{X}(t))|^{\bar{p}} \leq Ch^{\bar{p}}(k(h))^{2\bar{p}}.
    \end{split}
\end{align*}
The proof is complete.
\eproof

\section{Main results}

\begin{theorem}
    \label{theorem3-2}
    Let Assumptions \ref{ass1}, \ref{ass2} and \ref{ass4} hold, and let Assumption \ref{ass3} hold for any $q>2$, then for any $\bar{p}\in [2,p)$ and $\varepsilon\in(0,\frac{1}{4}]$, there exists a constant C such that for any $h \in (0,1]$  and $\l > 0$,
    \begin{align}
        \label{th311}
        \EE\left( \sup_{0\leq t\leq T}|Y(t)-X(t)|^{\bar{p}}\right)\leq h^{\min\{\gamma_f\bar{p},\gamma_g\bar{p},(1-2\varepsilon)\bar{p}\}}
    \end{align}
    and
    \begin{align}
        \label{th312}
        \EE\left( \sup_{0\leq t\leq T}|Y(t)-\bar{X}(t)|^{\bar{p}}\right)\leq Ch^{\min\{\gamma_f\bar{p},\gamma_g\bar{p},(1-2\varepsilon)\bar{p}\}}.
    \end{align}
\end{theorem}

\noindent
{\bf Proof.} Fix $\bar{p}\in [2,p)$ and $h\in(0,1]$ arbitrarily. Let $e(t)=Y(t)-X(t)$ for $t\geq 0$. For each integer $\ell> |Y(0)|$, define the stopping time
\begin{align}
    \theta_{\ell}=\inf\{t\geq0:|Y(t)|\vee|X(t)|\geq \ell\},
\end{align}
where we set $\inf\emptyset=\infty$ (as usual $\emptyset$ denotes the empty set). By the It\^o formula, we have that for any $0\leq t\leq T$,
\begin{align}
    \label{error}
        |e(t\wedge \theta_{\ell})|^{\bar{p}}&=\int^{t\wedge\theta_{\ell}}_{0}\bigg(\bar{p}|e(s)|^{\bar{p}-1}\left(f(s,Y(s))-f_{h}(\bar{\tau}(s),\bar{X}(s))\right)\nonumber\\
        &\quad +\frac{\bar{p}(\bar{p}-1)}{2}|e(s)|^{\bar{p}-2}\big|g(s,Y(s))-g_{h}(\bar{\tau}(s),\bar{X}(s))\nonumber\\
        &\quad
        -Lg_{h}(\bar{\tau}(s),\bar{X}(s))\Delta W(E_h(s))\big|^2 \bigg)dE(s)+M_{t\wedge \theta_{\ell}},
\end{align}
where
\begin{align*}
    M_{t\wedge \theta_{\ell}}:= &\int^{t\wedge\theta_{\ell}}_{0}\bar{p}|e(s)|^{\bar{p}-1}\big|g(s,Y(s))-g_{h}(\bar{\tau}(s),\bar{X}(s))\\
    &-Lg_{h}(\bar{\tau}(s),\bar{X}(s))\Delta W(E_h(s))\big|dW(E(s)).
\end{align*}
Note that the stochastic integral $(M_t)_{t\geq 0}$ is a local martingale with quadratic variation
\begin{align*}
    [M,M]_{t\wedge \theta_{\ell}}=&\int^{t\wedge\theta_{\ell}}_{0}\bar{p}^2|e(s)|^{2\bar{p}-2}\big|g(s,Y(s))-g_{h}(\bar{\tau}(s),\bar{X}(s))\\
    &-Lg_{h}(\bar{\tau}(s),\bar{X}(s))\Delta W(E_h(s))\big|^2 dE(s).
\end{align*}

For $0\leq s\leq t\wedge\theta_{\ell}$, we have
\begin{align*}
    \begin{split}
        &\bar{p}^2|e(s)|^{2\bar{p}-2}\big|g(s,Y(s))-g_{h}(\bar{\tau}(s),\bar{X}(s))-Lg_{h}(\bar{\tau}(s),\bar{X}(s))\Delta W(E_h(s))\big|^2\\
        =&\bar{p}^2|e(s)|^{\bar{p}}|e(s)|^{\bar{p}-2}\big|g(s,Y(s))-g_{h}(\bar{\tau}(s),\bar{X}(s))-Lg_{h}(\bar{\tau}(s),\bar{X}(s))\Delta W(E_h(s))\big|^2\\
        \leq& \bar{p}^2(\sup_{0\leq r\leq t\wedge\theta_{\ell}}|e(r)|^{\bar{p}})|e(s)|^{\bar{p}-2}\big|g(s,Y(s))-g_{h}(\bar{\tau}(s),\bar{X}(s))\\
        &-Lg_{h}(\bar{\tau}(s),\bar{X}(s))\Delta W(E_h(s))\big|^2.
    \end{split}
\end{align*}
Hence, using the inequality $(ab)^{1/2}\leq a/\l+\l b$ valid for any $a,b>0$ and $\l>0$, with $\l=2\bar{p}$ we have
\begin{align}
    \label{mt}
        &([M,M]_{t\wedge \theta_{\ell}})^{1/2}\nonumber\\
        \leq &\bar{p}\bigg(\sup_{0\leq r\leq t\wedge\theta_{\ell}}|e(r)|^{\bar{p}}\int^{t\wedge\theta_{\ell}}_{0}|e(s)|^{\bar{p}-2}\big|g(s,Y(s))-g_{h}(\bar{\tau}(s),\bar{X}(s))\nonumber\\
        \quad&
        -Lg_{h}(\bar{\tau}(s),\bar{X}(s))\Delta W(E_h(s))\big|^2dE(s)\bigg)^{\frac{1}{2}}\nonumber\\
        \leq & \bar{p}\bigg(\frac{\sup_{0\leq r\leq t\wedge\theta_{\ell}}|e(r)|^{\bar{p}}}{2\bar{p}}+2\bar{p}\int^{t\wedge\theta_{\ell}}_{0}|e(s)|^{\bar{p}-2}\big|g(s,Y(s))\nonumber\\
        \quad&
        -g_{h}(\bar{\tau}(s),\bar{X}(s))-Lg_{h}(\bar{\tau}(s),\bar{X}(s))\Delta W(E_h(s))\big|^2dE(s)\bigg)\nonumber\\
        \leq &\frac{1}{2}\sup_{0\leq r\leq t\wedge\theta_{\ell}}|e(r)|^{\bar{p}}+2\bar{p}^{2}\int^{t\wedge\theta_{\ell}}_{0}\bigg(|e(s)|^{\bar{p}-2}\big|g(s,Y(s))\nonumber\\
        \quad&
        -g_{h}(\bar{\tau}(s),{X}(s))+\tilde{R}_{g_h}(s,X(s),\bar{X}(s))\big|^2\bigg)dE(s),
\end{align}
where \eqref{the2_12} is used to get the last inequality.

We have expectation from \eqref{mt} 
\begin{align}
        \EE_W(M_{t\wedge \theta_{\ell}})=&\EE_W([M,M]_{t\wedge \theta_{\ell}})^{\frac{1}{2}}\nonumber\\
        =&\EE_W\bigg(\frac{1}{2}\sup_{0\leq r\leq t\wedge\theta_{\ell}}|e(r)|^{\bar{p}}+2\bar{p}^{2}\int^{t\wedge\theta_{\ell}}_{0}|e(s)|^{\bar{p}-2}\big|g(s,Y(s))\nonumber\\
        &-g_{h}(\bar{\tau}(s),{X}(s))+\tilde{R}_{g_h}(s,X(s),\bar{X}(s))\big|^2dE(s)\bigg).
\end{align}
Combing \eqref{error} and \eqref{mt} then we have
\begin{align*}
    \label{err226}
        &\EE_W\left(\sup_{0\leq t\leq T}|e(t\wedge \theta_{\ell})|^{\bar{p}}\right)\\
        \leq&\EE_W\bigg(\sup_{0\leq t\leq T}\int^{t\wedge\theta_{\ell}}_{0}\bar{p}|e(s)|^{\bar{p}-2}\bigg(|e(s)|^{\mathrm{T}}\left(f(s,Y(s))-f_{h}(\bar{\tau}(s),\bar{X}(s))\right)\\
        \quad& +\frac{\bar{p}-1}{2}\big|g(s,Y(s))-g_{h}(\bar{\tau}(s),\bar{X}(s))-Lg_{h}(\bar{\tau}(s),\bar{X}(s))\Delta W(E_h(s))\big|^2 \bigg)dE(s)\\
        \quad&
        +\frac{1}{2}\sup_{0\leq r\leq t\wedge\theta_{\ell}}|e(r)|^{\bar{p}}+2\bar{p}^{2}\sup_{0\leq t\leq T}\int^{t\wedge\theta_{\ell}}_{0}|e(s)|^{\bar{p}-2}\big|g(s,Y(s))-g_{h}(\bar{\tau}(s),X(s))\\
        \quad&
        +\tilde{R}_{g_h}(s,X(s),\bar{X}(s))\big|^2dE(s)\bigg)\\
        \leq&\EE_W\bigg(\sup_{0\leq t\leq T}\int^{t\wedge\theta_{\ell}}_{0}\bar{p}|e(s)|^{\bar{p}-2}\bigg(|e(s)|^{\mathrm{T}}\left(f(s,Y(s))-f_{h}(\bar{\tau}(s),\bar{X}(s))\right)\\
        \quad& +\frac{\bar{p}-1}{2}\big|g(s,Y(s))-g_{h}(\bar{\tau}(s),X(s))+\tilde{R}_{g_h}(s,X(s),\bar{X}(s))\big|^2\bigg)dE(s)\\
        \quad&
        +\frac{1}{2}\sup_{0\leq r\leq t\wedge\theta_{\ell}}|e(r)|^{\bar{p}}+2\bar{p}^{2}\sup_{0\leq t\leq T}\int^{t\wedge\theta_{\ell}}_{0}|e(s)|^{\bar{p}-2}\big|g(s,Y(s))-g_{h}(\bar{\tau}(s),X(s))\\
        \quad&
        +\tilde{R}_{g_h}(s,X(s),\bar{X}(s))\big|^2dE(s)\bigg).
\end{align*}
Where the second term uses \eqref{the2_12}, then, by organizing the above equations, we obtained,
\begin{align}
        &\EE_W\left(\sup_{0\leq t\leq T}|e(t)\wedge \theta_{\ell})|^{\bar{p}}\right)\nonumber\\
        \leq&\EE_W\bigg(\sup_{0\leq t\leq T}\int^{t\wedge\theta_{\ell}}_{0}\bar{p}|e(s)|^{\bar{p}-2}\bigg(|e(s|^{\mathrm{T}}\left(f(s,Y(s))-f_{h}(\bar{\tau}(s),\bar{X}(s))\right)\nonumber\\
        \quad& +(\bar{p}-1)|g(s,Y(s))-g_{h}(\bar{\tau}(s),X(s))|^{2}\nonumber\\
        \quad &+(\bar{p}-1)
        |\tilde{R}_{g_h}(s,X(s),\bar{X}(s))|^2\bigg)dE(s)+\frac{1}{2}\sup_{0\leq r\leq t\wedge\theta_{\ell}}|e(r)|^{\bar{p}}\nonumber\\
        \quad &+\bar{p}|e(s)|^{\bar{p}-2}\sup_{0\leq t\leq T}\int^{t\wedge\theta_{\ell}}_{0}\bigg(4\bar{p}|g(s,Y(s))-g_{h}(\bar{\tau}(s),X(s))|^{2}\nonumber\\
        \quad&
        +4\bar{p}|\tilde{R}_{g_h}(s,X(s),\bar{X}(s))|^2\bigg)dE(s)\bigg)\nonumber\\
        \leq&\EE_W\bigg(\sup_{0\leq t\leq T}\int^{t\wedge\theta_{\ell}}_{0}\bar{p}|e(s)|^{\bar{p}-2}\bigg(|e(s)|^{\mathrm{T}}\left(f(s,Y(s))-f_{h}(\bar{\tau}(s),\bar{X}(s))\right)\nonumber\\
        \quad& +(5\bar{p}-1)|g(s,Y(s))-g_{h}(\bar{\tau}(s),X(s))|^{2}\bigg)dE(s)+\frac{1}{2}\sup_{0\leq r\leq t\wedge\theta_{\ell}}|e(r)|^{\bar{p}}\nonumber\\
        \quad&
        +\bar{p}(5\bar{p}-1)\sup_{0\leq t\leq T}\int^{t\wedge\theta_{\ell}}_{0}|e(s)|^{\bar{p}-2} |\tilde{R}_{g_h}(s,X(s),\bar{X}(s))|^2dE(s)\bigg).   
\end{align}
For the last two items uses a basic inequality $(a+b)^{2}\leq2(a^{2}+b^{2})$ and then merge.

Next, Let's organize the equation and use the Young inequality $(a+b)^{2}\leq (1+\varepsilon)a^{2}+(1+1/\varepsilon)b^{2}$ for any $a,b\geq 0$, $ \varepsilon>0$, we choose $ \varepsilon =(5p-5\bar{p})/(5\bar{p}-1)$ in the second term,  we can get from \eqref{err226}
\begin{align}
    \label{err2}
        &\EE_W\left(\sup_{0\leq t\leq T}|e(t\wedge \theta_{\ell})|^{\bar{p}}\right)\nonumber\\
        \leq&\EE_W\bigg(\sup_{0\leq t\leq T}\int^{t\wedge\theta_{\ell}}_{0}\bar{p}|e(s)|^{\bar{p}-2}\bigg(|e(s)|^{\mathrm{T}}\left(f(s,Y(s))-f_{h}(\bar{\tau}(s),\bar{X}(s))\right)\nonumber\\
        \quad& +(5\bar{p}-1)\big|g(s,Y(s))-g(s,X(s))+g(s,X(s))\nonumber\\
        \quad&
        -g_{h}(\bar{\tau}(s),X(s))\big|^{2}\bigg)dE(s)+\frac{1}{2}\sup_{0\leq r\leq t\wedge\theta_{\ell}}|e(r)|^{\bar{p}}\nonumber\\
        \quad&+(5{\bar{p}}^{2}-\bar{p})\sup_{0\leq t\leq T}\int^{t\wedge\theta_{\ell}}_{0}|e(s)|^{\bar{p}-2}|\tilde{R}_{g_h}(s,X(s),\bar{X}(s))|^2dE(s)\bigg)\nonumber\\
        \leq&\EE_W\bigg(\sup_{0\leq t\leq T}\int^{t\wedge\theta_{\ell}}_{0}\bar{p}|e(s)|^{\bar{p}-2}\bigg(|e(s)|^{\mathrm{T}}\left(f(s,Y(s))-f_{h}(\bar{\tau}(s),\bar{X}(s))\right)\nonumber\\
        \quad& +(5\bar{p}-1)\bigg((1+\dfrac{5p-5\bar{p}}{5\bar{p}-1})|g(s,Y(s))-g(s,X(s))|^{2}\nonumber\\
        \quad&
        +(1+\dfrac{5\bar{p}-1}{5p-5\bar{p}})|g(s,X(s))-g_{h}(\bar{\tau}(s),X(s))|^{2}\bigg)dE(s)\nonumber\\
        \quad&
        +\frac{1}{2}\sup_{0\leq r\leq t\wedge\theta_{\ell}}|e(r)|^{\bar{p}}+(5{\bar{p}}^{2}-\bar{p})\sup_{0\leq t\leq T}\int^{t\wedge\theta_{\ell}}_{0}|e(s)|^{\bar{p}-2}\nonumber\\
        \quad& \times|\tilde{R}_{g_h}(s,X(s),\bar{X}(s))|^2dE(s)\bigg)\nonumber\\
        \leq&\EE_W\bigg(\sup_{0\leq t\leq T}\int^{t\wedge\theta_{\ell}}_{0}\bar{p}|e(s)|^{\bar{p}-2}\bigg(|e(s)|^{\mathrm{T}}\left(f(s,Y(s))-f_{h}(\bar{\tau}(s),\bar{X}(s))\right)\nonumber\\
        \quad& +(5p-1)|g(s,Y(s))-g(s,X(s))|^{2}
        +\dfrac{5p-1}{5p-5\bar{p}}|g(s,X(s))\nonumber\\
        \quad&
        -g_{h}(\bar{\tau}(s),X(s))|^{2}\bigg)dE(s)+\frac{1}{2}\sup_{0\leq r\leq t\wedge\theta_{\ell}}|e(r)|^{\bar{p}}\nonumber\\
        \quad&+(5{\bar{p}}^{2}-\bar{p})\sup_{0\leq t\leq T}\int^{t\wedge\theta_{\ell}}_{0}|e(s)|^{\bar{p}-2}|\tilde{R}_{g_h}(s,X(s),\bar{X}(s))|^2dE(s)\bigg).
\end{align}
Using the basic properties of inequalities, we get from\eqref{err2} that
\begin{align}
    \label{sup}
        &\EE_W\left(\sup_{0\leq t\leq T}|e( t\wedge\theta_{\ell})|^{\bar{p}}\right)\nonumber\\
        \leq& \frac{1}{2}\sup_{0\leq r\leq t\wedge\theta_{\ell}}|e(r)|^{\bar{p}}
        +\EE_W\sup_{0\leq t\leq T}\int^{t\wedge\theta_{\ell}}_{0}\bar{p}|e(s)|^{\bar{p}-2}\bigg(e^{\mathrm{T}}(s)\big(f(s,Y(s))\nonumber\\
        &-f(s,X(s))\big)
        +(5p-1)|g(s,Y(s))-g(s,X(s))|^2\bigg)dE(s)\nonumber\\
        &+\EE_W\sup_{0\leq t\leq T}\int^{t\wedge\theta_{\ell}}_{0}\bar{p}|e(s)|^{\bar{p}-2}\bigg(e^{\mathrm{T}}(s)\big(f(s,X(s))-f_{h}(\bar{\tau}(s),X(s))\big)\nonumber\\
        \quad &+\dfrac{5p-1}{5p-5\bar{p}}|g(s,X(s))-g_{h}(\bar{\tau}(s),X(s))|^2\bigg)dE(s)\nonumber\\
        &+\EE_W\sup_{0\leq t\leq T}\int^{t\wedge\theta_{\ell}}_{0}(5\bar{p}^{2}-p)|e(s)|^{\bar{p}-2}|\tilde{R}_{g_h}(s,X(s),\bar{X}(s))|^2dE(s)\nonumber\\
        \leq& \frac{1}{2}\sup_{0\leq r\leq t\wedge\theta_{\ell}}|e(r)|^{\bar{p}}+[J_1]+[J_2]+[J_3],
\end{align}
where
\begin{align*}
        J_1&:=\EE_W \bigg(\sup_{0\leq t\leq T}\int^{t\wedge\theta_{\ell}}_{0}\bar{p}|e(s)|^{\bar{p}-2}\bigg(e^{\mathrm{T}}(s)\big(f(s,Y(s))-f(s,X(s))\big)\\
        &\quad+(5p-1)|g(s,Y(s))-g(s,X(s))|^2\bigg)dE(s)\bigg),
\end{align*}
\begin{align*}
        J_2&:=\EE_W\bigg(\sup_{0\leq t\leq T}\int^{t\wedge\theta_{\ell}}_{0}\bar{p}|e(s)|^{\bar{p}-2}\bigg(e^{\mathrm{T}}(s)\left(f(s,X(s))-f_{h}(\bar{\tau}(s),\bar{X}(s))\right)\\
        &\quad+\frac{5p-1}{5p-5\bar{p}}|g(s,X(s))-g_{h}(\bar{\tau}(s),X(s))|^2\bigg)dE(s)\bigg),
\end{align*}
\begin{align*}
    \begin{split}
        J_3&:=\EE_W\bigg(\sup_{0\leq t\leq T}\int^{t\wedge\theta_{\ell}}_{0}(5\bar{p}^{2}-p)|e(s)|^{\bar{p}-2}|\tilde{R}_{g_h}(s,X(s),\bar{X}(s))|^2dE(s)\bigg).
    \end{split}
\end{align*}
By Assumption \ref{ass2}, we have
\begin{align}
    \label{J1}
    J_1\leq H_1\int^{T}_{0}\EE_W|e(s)|^{\bar{p}}dE(s),
\end{align}
where $H_1=\bar{p}K$. Next, handling the $J_{2}$
\begin{align}
    \label{lemm330}
        J_2&= \EE_W\bigg(\sup_{0\leq t\leq T}\int^{t\wedge\theta_{\ell}}_{0}\bar{p}|e(s)|^{\bar{p}-2}\bigg(e^{T}(s)\big(f(s,X(s))-f_h(\bar{\tau}(s),\bar{X}(s))\big)\nonumber\\
        &\quad+\dfrac{5p-1}{5p-5\bar{p}}|g(s,X(s))-g_h(\bar{\tau}(s),X(s))|^2\bigg)dE(s)\bigg)\nonumber\\
        &\leq \EE_W\bigg(\sup_{0\leq t\leq T}\int^{t\wedge\theta_{\ell}}_{0}\bar{p}|e(s)|^{\bar{p}-2}\bigg(e^{T}(s)\big(f(s,X(s))-f(\bar{\tau}(s),X(s))\big)\nonumber\\
        &\quad+\dfrac{5p-1}{5p-5\bar{p}}|g(s,X(s))-g(\bar{\tau}(s),X(s))|^2\bigg)dE(s)\nonumber\\
        &\quad+\sup_{0\leq t\leq T}\int^{t\wedge\theta_{\ell}}_{0}\bar{p}|e(s)|^{\bar{p}-2}\bigg(e^{T}(s)\big(f(\bar{\tau}(s),X(s))-f_h(\bar{\tau}(s),\bar{X}(s))\big)\nonumber\\
        &\quad+\dfrac{5p-1}{5p-5\bar{p}}|g(\bar{\tau}(s),X(s))-g_h(\bar{\tau}(s),X(s))|^2\bigg)dE(s)\bigg)\nonumber\\
        &\leq J_{21}+J_{22},
\end{align}
where 
\begin{align*}
        J_{21}&= \EE_W\bigg(\sup_{0\leq t\leq T}\int^{t\wedge\theta_{\ell}}_{0}\bar{p}|e(s)|^{\bar{p}-2}\bigg(e^{T}(s)\big(f(s,X(s))-f(\bar{\tau}(s),X(s))\big)\\
        &\quad+\dfrac{5p-1}{5p-5\bar{p}}|g(s,X(s))-g(\bar{\tau}(s),X(s))|^2\bigg)dE(s)\bigg),
\end{align*}
\begin{align*}
        J_{22}&= \EE_W\bigg(\sup_{0\leq t\leq T}\int^{t\wedge\theta_{\ell}}_{0}\bar{p}|e(s)|^{\bar{p}-2}\bigg(e^{T}(s)\big(f(\bar{\tau}(s),X(s))-f_h(\bar{\tau}(s),\bar{X}(s))\big)\\
        &\quad+\dfrac{5p-1}{5p-5\bar{p}}|g(\bar{\tau}(s),X(s))-g_h(\bar{\tau}(s),X(s))|^2\bigg)dE(s)\bigg).
\end{align*}
Using Assumption \ref{ass4}, basic inequality and the Young inequality, for any $0\leq t\leq t\wedge\theta_{\ell}\leq T $,
\begin{align*}
    a^{p-2}b\leq \frac{p-2}{p}a^p + \frac{2}{p}b^{p/2},\quad \forall a,b \geq 0.
\end{align*}
We can derive 
\begin{align}
    \label{lemm331}
        J_{21}&\leq
        \EE_W\bigg(\sup_{0\leq t\leq T}\int^{t\wedge\theta_{\ell}}_{0}\bar{p}|e(s)|^{\bar{p}-2}\bigg(\frac{1}{2}|e(s)|^2+\frac{1}{2}|f(s,X(s))-f(\bar{\tau}(s),X(s))|^2\nonumber\\
        &\quad+\dfrac{5p-1}{5p-5\bar{p}}|g(s,X(s))-g(\bar{\tau}(s),X(s))|^2\bigg)dE(s)\bigg)\nonumber\\
        &\leq
        C\bigg(\EE_W\sup_{0\leq t\leq T}\bigg(\int^{t\wedge\theta_{\ell}}_{0}|e(s)|^{\bar{p}}dE(s)+ \int^{t\wedge\theta_{\ell}}_{0}|f(s,X(s))-f(\bar{\tau}(s),X(s))|^{\bar{p}}dE(s)\nonumber\\
        &\quad +\int^{t\wedge\theta_{\ell}}_{0}|g(s,X(s))-g(\bar{\tau}(s),X(s))|^{\bar{p}}\bigg)dE(s)\bigg)\nonumber\\
        &\leq
        C\bigg(\EE_W\int^{T}_{0}|e(s)|^{\bar{p}}dE(s)
        +\EE_W\int^{T}_{0}H_1^{\bar{p}}(1+|X(s)|^{(1+\a)\bar{p}})h ^{\gamma_f\bar{p}}dE(s)\nonumber\\
        &\quad  +\EE_W\int^{T}_{0}H_2^{\bar{p}}(1+|X(s)|^{(1+\a)\bar{p}})h^{\gamma_g\bar{p}})dE(s)\bigg)\nonumber\\
        &\leq
        C\bigg(\EE_W\int^{T}_{0}|e(s)|^{\bar{p}}dE(s)
        +h^{\gamma_f\bar{p}}E(T)+h^{\gamma_g\bar{p}}E(T)\bigg).
\end{align}
Where lemma \ref{lemma3} is used to get the last inequality.\\
We use the basic properties of inequalities to handle the $J_{22}$ item,
\begin{align}
    \label{lem332}
        J_{22}&\leq  \EE_W\bigg(\sup_{0\leq t\leq T}\int^{t\wedge\theta_{\ell}}_{0}\bar{p}|e(s)|^{\bar{p}-2}\bigg(e^{T}(s)\big(f(\bar{\tau}(s),X(s))-f(\bar{\tau}(s),\bar{X}(s))\big)\bigg)dE(s)\nonumber\\
        &\quad+
        \sup_{0\leq t\leq T}\int^{t\wedge\theta_{\ell}}_{0}\bar{p}|e(s)|^{\bar{p}-2}\bigg(e^{T}(s)\big(f(\bar{\tau}(s),\bar{X}(s))-f_{h}(\bar{\tau}(s),\bar{X}(s))\big)\nonumber\\
        &\quad+\dfrac{5p-1}{5p-5\bar{p}}|g(\bar{\tau}(s),X(s))-g_h(\bar{\tau}(s),X(s))|^2\bigg)dE(s)\bigg)\nonumber\\
        &\leq I_1+ I_2.
\end{align}
We can derive from the \eqref{le210} and Young inequality
\begin{align}
        \label{333}
        I_{1}&=\EE_W\bigg(\sup_{0\leq t\leq T}\int^{t\wedge\theta_{\ell}}_{0}\bar{p}|e(s)|^{\bar{p}-2}\bigg(e^{T}(s)\big(f(\bar{\tau}(s),X(s))-f(\bar{\tau}(s),\bar{X}(s))\big)\bigg)dE(s)\bigg)\nonumber\\
        &\leq  \EE_W\bigg(\sup_{0\leq t\leq T}\int^{t\wedge\theta_{\ell}}_{0}\bar{p}|e(s)|^{\bar{p}-2}\bigg(e^{T}(s)\big(f^{'}(\bar{\tau}(s),x)|_{x=\bar{X}(s)}\int^{s}_{0}g_h(\bar{\tau}(s_{1}),\bar{X}(s_{1}))\nonumber\\
        &\quad \times dW(E(s_{1}))+\tilde{R}_{f}(s,X(s),\bar{X}(s))\big)\bigg)dE(s)\bigg)\nonumber\\
        &\leq 
        H_{21}\EE_W\bigg(\sup_{0\leq t\leq T}\int^{t\wedge\theta_{\ell}}_{0}\bigg(|e(s)|^{\bar{p}}+\big|e(s)^{T}(f^{'}(\bar{\tau}(s),x)|_{x=\bar{X}(s)}\int^{s}_{0}g_h(\bar{\tau}(s_{1}),\bar{X}(s_{1}))\nonumber\\
        &\quad \times dW(E(s_{1}))\big|^{\frac{\bar{p}}{2}}+|e(s)^{T}\tilde{R}_{f}(s,X(s),\bar{X}(s))|^{{\frac{\bar{p}}{2}}}\bigg)dE(s)\bigg)\nonumber\\
        &\leq  H_{21}\bigg(\EE_W\sup_{0\leq t\leq T}\int^{t\wedge\theta_{\ell}}_{0}\bigg(|e(s)|^{\bar{p}}dE(s)+\big|e(s)^{T}(f^{'}(\bar{\tau}(s),x)|_{x=\bar{X}(s)}\int^{s}_{0}g_h(\bar{\tau}(s_{1}),\bar{X}(s_{1}))\nonumber\\
        &\quad \times dW(E(s_{1}))\big|^{\frac{\bar{p}}{2}}dE(s)+|\tilde{R}_{f}(s,X(s),\bar{X}(s))|^{\bar{p}}dE(s)\bigg).   
\end{align}
Apply a similar approach, used for (3.35) in \cite{Wang2013} and combing \eqref{333} and lemma \ref{lemma5}, we obtain
\begin{align}
    \label{335}
        I_{1}&\leq H_{21}\bigg(\EE_{W}\int^{T}_{0}|e(s)|^{\bar{p}}dE(s)+\EE_{W}\int^{T}_{0}|\tilde{R}_{f}(s,X(s),\bar{X}(s))|^{\bar{p}}dE(s)+h^{\bar{p}}\bigg)\nonumber\\
        &\leq 
        H_{21}\bigg(\EE_{W}\int^{T}_{0}|e(s)|^{\bar{p}}dE(s)+\int^{T}_{0}\EE_{W}|\tilde{R}_{f}(s,X(s),\bar{X}(s))|^{\bar{P}}dE(s)+h^{\bar{p}}\bigg)\nonumber\\
        &\leq 
        H_{21}\bigg(\EE_{W}\int^{T}_{0}|e(s)^{\bar{p}}dE(s))+h^{\bar{p}}(k(h))^{2\bar{p}}+h^{\bar{p}}\bigg).
\end{align} 
Applying the Young inequality, Assumption \ref{ass1} and H{\"o}lder inequality, we can show that
\begin{align*}
    I_{2}&=\EE_{W}\bigg(\sup_{0\leq t\leq T}\int^{t\wedge\theta_{\ell}}_{0}\bar{p}|e(s)|^{\bar{p}-2}\bigg(e^{T}(s)\big(f(\bar{\tau}(s),\bar{X}(s))-f_h(\bar{\tau}(s),\bar{X}(s))\big)\\
    &\quad+\dfrac{5p-1}{5p-5\bar{p}}|g(\bar{\tau}(s),X(s))-g_h(\bar{\tau}(s),X(s))|^2\bigg)dE(s)\bigg)\\
    &\leq 
    H_{22}\bigg(\EE_{W}\sup_{0\leq t\leq T}\int^{t\wedge\theta_{\ell}}_{0}|e(s)|^{\bar{p}}dE(s)+\EE_{W}\sup_{0\leq t\leq T}\int^{t\wedge\theta_{\ell}}_{0}|f(\bar{\tau}(s),\bar{X}(s))-f_h(\bar{\tau}(s),\bar{X}(s))|^{\bar{p}}\\
    &\quad+|g(\bar{\tau}(s),X(s))-g_h(\bar{\tau}(s),X(s))|^{\bar{p}}dE(s)\bigg)\\
    &\leq 
    H_{22}\bigg(\EE_{W}\sup_{0\leq t\leq T}\int^{t\wedge\theta_{\ell}}_{0}|e(s)|^{\bar{p}}dE(s)+\EE_{W}\sup_{0\leq t\leq T}\int^{t\wedge\theta_{\ell}}_{0}\big(1+|\bar{X}(s)|^{\a\bar{p}}\\
    &\quad +\big||\bar{X}(s)|\wedge \mu ^{-1}(k(h))\big|^{\a\bar{p}}\big) \bigg|\bar{X}(s)-\big(|\bar{X}(s)|\wedge \mu ^{-1}(k(h))\big)\frac{\bar{X}(s)}{|\bar{X}(s)|}\bigg|^{\bar{p}}dE(s)\\
    &\quad + \EE_{W}\sup_{0\leq t\leq T}\int^{t\wedge\theta_{\ell}}_{0}(1+|X(s)|^{\a\bar{p}}+\big||X(s)|\wedge \mu ^{-1}(k(h))\big|^{\a\bar{p}})\\
    &\quad \times \bigg|X(s)-\big(|X(s)|\wedge \mu ^{-1}(k(h))\big)\frac{X(s)}{|X(s)|}\bigg|^{\bar{p}}dE(s) \bigg)\\
    &\leq 
    H_{22}\bigg(\EE_{W}\int^{T}_{0}|e(s)|^{\bar{p}}dE(s)+\int^{T}_{0}\bigg(\EE_{W}\bigg[1+|\bar{X}(s)|^{q}+\big||\bar{X}(s)|\wedge \mu ^{-1}(k(h))\big|^{q}\bigg]\bigg)^{\frac{\a\bar{p}}{q}}\\
    &\quad \times \bigg[\EE_{W}\big|\bar{X}(s)-\big(|\bar{X}(s)|\wedge \mu ^{-1}(k(h))\big)\frac{\bar{X}(s)}{|\bar{X}(s)|}\big|^{\frac{q\bar{p}}{q-\a\bar{p}}}\bigg]^{\frac{q-\a\bar{p}}{q}}dE(s)\\
    &\quad + \int^{T}_{0}\bigg(\EE_{W}\bigg[1+|X(s)|^{q}+\big||X(s)|\wedge \mu ^{-1}(k(h))\big|^{q}\bigg]\bigg)^{\frac{\a\bar{p}}{q}}\\
    &\quad \times\bigg[\EE_{W} \big|X(s)-\big(|X(s)|\wedge \mu ^{-1}(k(h))\big)\frac{X(s)}{|X(s)|}\big|^{\frac{q\bar{p}}{q-\a\bar{p}}}\bigg]^{\frac{q-\a\bar{p}}{q}}dE(s) \bigg),
\end{align*}
where the lemma \ref{lemma3} are used above, also used in the follwing last inequality, using the H{\"o}lder inequality and chebyshev inequality $ \PP(|x|\geqslant a) \leqslant a^{-q} \EE|x|^{q}$, if $a>0,q>0$ , we can obtain\\
\begin{align}
        \label{336}
        I_{2}&\leq 
        H_{22}\bigg(\EE_{W}\int^{T}_{0}|e(s)|^{\bar{p}}dE(s)\nonumber\\
        &\quad+\int^{T}_{0}\bigg(\EE_{W}\big|I\left\{|\bar{X}(s)|>\mu^{-1}(k(h))\right\}|\bar{X}(s)|^{\frac{q\bar{p}}{q-\a\bar{p}}}\big|\bigg)^{\frac{q-\a\bar{p}}{q}}dE(s)\nonumber\\
        &\quad+\int^{T}_{0}\bigg(\EE_{W}\big|I\left\{|X(s)|>\mu^{-1}(k(h))\right\}|X(s)|^{\frac{q\bar{p}}{q-\a\bar{p}}}\big|\bigg)^{\frac{q-\a\bar{p}}{q}}dE(s)\bigg)\nonumber\\
        &\leq 
        H_{22}\bigg(\EE_{W}\int^{T}_{0}|e(s)|^{\bar{p}}dE(s)\nonumber\\
        &\quad+\int^{T}_{0}\bigg(\big[P\left\{|\bar{X}(s)|>\mu^{-1}(k(h))\right\}\big]^{\frac{q-\a\bar{p}-\bar{p}}{q-\a\bar{p}}} \big[\EE|\bar{X}(s)|^{q}\big]^{\frac{\bar{p}}{q-\a\bar{p}}}\bigg)^{\frac{q-\a\bar{p}}{q}}dE(s)\nonumber\\
        &\quad+\int^{T}_{0}\bigg(\big[P\left\{|X(s)|>\mu^{-1}(k(h))\right\}\big]^{\frac{q-\a\bar{p}-\bar{p}}{q-\a\bar{p}}}[\EE|X(s)|^{q}]^{\frac{\bar{p}}{q-\a\bar{p}}}\bigg)^{\frac{q-\a\bar{p}}{q}}dE(s)\bigg)\nonumber\\
        &\leq 
        H_{22}\bigg(\EE_{W}\int^{T}_{0}|e(s)|^{\bar{p}}dE(s)+\int^{T}_{0}\bigg(\frac{\EE_{W}|\bar{X}(s)|^{q}}{|\mu^{-1}(k(h))|^{q}}\bigg)^{\frac{q-\a\bar{p}-\bar{p}}{q}}dE(s)\nonumber\\
        &\quad+\int^{T}_{0}\bigg(\frac{\EE_{W}|X(s)|^{q}}{|\mu^{-1}(k(h))|^{q}}\bigg)^{\frac{q-\a\bar{p}-\bar{p}}{q}}dE(s)\bigg)\nonumber\\
        &\leq 
        H_{22}\bigg(\EE_{W}\int^{T}_{0}|e(s)|^{\bar{p}}dE(s)+\big(\mu^{-1}(k(h))\big)^{(\a+1)\bar{p}-q}\bigg).
\end{align}
Substituing \eqref{335} and \eqref{336} into \eqref{lem332} gives
\begin{align}
        \label{lemm337}
        J_{22}\leq  H_{22}\bigg(\EE_{W}\int^{T}_{0}|e(s)|^{\bar{p}}dE(s)+\big(\mu^{-1}(k(h))\big)^{(\a+1)\bar{p}-q}+h^{\bar{p}}(k(h))^{2\bar{p}}+h^{\bar{p}}\bigg).
\end{align} 
Applying Young inequality and lemma \ref{lemma5}, we derive that
\begin{align}
    \label{J3}
        J_3&= \EE_{W}\sup_{0\leq t\leq T}\int^{t\wedge\theta_{\ell}}_{0}(5\bar{p}^{2}-\bar{p})|e(s)|^{\bar{p}-2}|\tilde{R}_{g_h}(s,X(s),\bar{X}(s))|^2dE(s)\nonumber\\
        &\leq 
        H_{3}\EE_{W}\sup_{0\leq t\leq T}\int^{t\wedge\theta_{\ell}}_{0}(|e(s)|^{\bar{p}}+|\tilde{R}_{g_h}(s,X(s),\bar{X}(s))|^{\bar{p}})dE(s)\nonumber\\
        &\leq 
        H_{3}\bigg(\EE_{W}\int^{T}_{0}|e(s)|^{\bar{p}}dE(s)+\int^{T}_{0}\EE_{W}|\tilde{R}_{g_h}(s,X(s),\bar{X}(s))|^{\bar{p}}dE(s)\bigg)\nonumber\\
        &\leq 
        H_{3}\bigg(\EE_{W}\int^{T}_{0}|e(s)|^{\bar{p}}dE(s)+h^{\bar{p}}(k(h))^{2\bar{p}}\bigg).
\end{align}
Where $H_{21}, H_{22}, H_{3}$ and following C are generic constants indendent of $h$ that may change from line to line, combining \eqref{sup}, \eqref{J1}, \eqref{lemm330}, \eqref{lemm331}, \eqref{lemm337}, \eqref{J3},
we can get the original formual is
\begin{align*}
        \EE_W\left(\sup_{0\leq t\leq T}|e( t\wedge\theta_{\ell})|^{\bar{p}}\right)\leq& \frac{1}{2}\sup_{0\leq r\leq t\wedge\theta_{\ell}}|e(r)|^{\bar{p}}+[J_1]+[J_2]+[J_3]\\
        \leq&
        2([J_1]+[J_2]+[J_3])\\
        \leq  &
        C\bigg(\EE_{W}\int^{T}_{0}|e(s)|^{\bar{p}}dE(s)+h^{\gamma_f\bar{p}}+h^{\gamma_g\bar{p}}\\
        &\quad+h^{\bar{p}}(k(h))^{2\bar{p}}+h^{\bar{p}}+\big(\mu^{-1}(k(h))\big)^{(\a+1)\bar{p}-q}\bigg)\\
        \leq &
        C\bigg(\int^{T}_{0}\EE_{W}\sup_{0\leq u\leq s}|e(u\wedge\theta_{\ell} )|^{\bar{p}}dE(s)+h^{\gamma_f\bar{p}}+h^{\gamma_g\bar{p}}\\
        &\quad+h^{\bar{p}}(k(h))^{2\bar{p}}+\big(\mu^{-1}(k(h))\big)^{(\a+1)\bar{p}-q}\bigg).
\end{align*}
An application of the Gronwall inequality yields that
\begin{align*}
        \EE_W(\sup_{0\leq t\leq T}|e( t\wedge\theta_{\ell})|^{\bar{p}})
        &\leq 
        C\bigg(h^{\gamma_f\bar{p}}+h^{\gamma_g\bar{p}}
        +h^{\bar{p}}(k(h))^{2\bar{p}}+(\mu^{-1}(k(h)))^{(\a+1)\bar{p}-q}\bigg)e^{\l E(T)},
\end{align*}
therefore thanks to the Fatou lemma, the assertion is proved by letting $n\to \infty$.
\begin{align*}
        \EE_W(\sup_{0\leq t\leq T}|e(t)|^{\bar{p}})
        &\leq 
        C\bigg(h^{\gamma_f\bar{p}}+h^{\gamma_g\bar{p}}
        +h^{\bar{p}}(k(h))^{2\bar{p}}+\big(\mu^{-1}(k(h))\big)^{(\a+1)\bar{p}-q}\bigg)e^{\l E(T)},
\end{align*}
where $C$ is independent from $h$ and $\l > 0$.
\par
Taking $\EE_D$ on both sides gives 
\begin{align}
        \label{th321}
        &\EE \bigg( \sup_{0\leq t\leq T}|Y(t)-X(t)|^{\bar{p}} \bigg) \nonumber\\
        \leq& C\bigg(h^{\gamma_f\bar{p}}+h^{\gamma_g\bar{p}}+h^{\bar{p}}(\k(h))^{2\bar{p}}+(\m^{-1}(\k(h)))^{(\a+1)\bar{p}-q}\bigg).
    \end{align}
Lemma \ref{lemma2} together with \eqref{th321}  indicates
\begin{align}
        \label{th322}
        &\EE\bigg(\sup_{0\leq t\leq T}|Y(t)-\bar{X}(t)|^{\bar{p}}\bigg)\nonumber\\
        \leq& C\bigg(h^{\gamma_f\bar{p}}+h^{\gamma_g\bar{p}}+h^{\bar{p}}(\k(h))^{2\bar{p}}+(\m^{-1}(\k(h)))^{(\a+1)\bar{p}-q}\bigg).
    \end{align}
 At last, by properly choosing $\mu^{-1}(\cdot)$ and $\kappa(\cdot)$, the required assertions are obtained

\eproof

\section{Numerical examples}
In this section, we give two numerical examples. 
\begin{expl}\label{ex}
    Consider a one-dimensional time-changed SDE
    \begin{align}\label{ex1}
        \left\{
        \begin{array}{lr}
            dY(t)=\left([t(1-t)]^{\frac{1}{4}}Y(t)-Y^5(t)\right)dE(t)+\left([t(1-t)]Y^2(t)\right)dW(E(t)),&\\
            Y(0)=1,&
        \end{array}
        \right.
    \end{align}
\end{expl}
with  $T=1$, the drift and diffusion coefficients are $f(y)=[t(1-t)]^{\frac{1}{4}}y-y^5$ and $g(y)=[t(1-t)]y^2$, respectively. Clearly, both of them have continuous second-order derivatives and it is not hard to verify that Assumption \ref{ass1} and Assumption \ref{ass5} are satisfied with $\alpha =4$.
\par
For any $p>2$, we can see
\begin{align*}
        &(x-y)^{\mathrm{T}}(f(t,x)-f(t,y))+(5p-1)|g(t,x)-g(t,y)|^2\\
        =&
        (x-y)^T\bigg([t(1-t)]^{\frac{1}{4}}(x-y)-(x^5-y^5)\bigg)+(5p-1)\big|[t(1-t)](x^2-y^2)\big|^{2}\\
        \leq&
        (x-y)^2\bigg([t(1-t)]^{\frac{1}{4}}-(x^4+x^{3}y+x^{2}y^{2}+xy^{3}+y^4)+(5p-1)[t(1-t)]^{2}(x+y)^{2}\bigg).
\end{align*}
But 
\begin{align*}
        -(x^{3}y+xy^{3})=-xy(x^{2}+y^{2})\leq 0.5(x^{2}+y^{2})^{2}=0.5(x^{4}+y^{4})+x^{2}y^{2}.
\end{align*}
Hence
\begin{align*}
        &(x-y)^{\mathrm{T}}(f(t,x)-f(t,y))+(5p-1)|g(t,x)-g(t,y)|^2\\
        \leq&(x-y)^2\bigg([t(1-t)]^{\frac{1}{4}}-0.5(x^4+y^4)+2(5p-1)[t(1-t)]^{2}(x^{2}+y^{2})\bigg)\\
        \leq& K(x-y)^{2},
\end{align*}
where the Young inequality is used. Note that the last inequality is due to the fact that polynomials with the negative coefficients for the highest order term can always be bounded from above. This indicates that Assumption \ref{ass2} holds.

In the similar manner, for any $q>2$ and any $t\in [0,1]$, we have
\begin{align*}
        & x^{\mathrm{T}}f(t,x)+(5q-1)|g(t,x)|^2\\
        =& [t(1-t)]^{\frac{1}{4}}x-x^5+(5q-1)[t(1-t)]^{2}x^4\\
        \leq & K_{1}(1+|x|^2),
\end{align*}
which means that Assumption \ref{ass3} is satisfied.

Using the mean theorem for the temporal variable, Assumption \ref{ass4} are satisfied with $\gamma_f=\frac{1}{4}$, $\gamma_g=1$. According to Theorem \ref{theorem3-2}, we know that
\begin{align*}
    \EE \left( \sup_{0\leq t\leq T}|Y(t)-X(t)|^{\bar{p}} \right) \leq C\bigg(h^{\frac{\bar{p}}{4}}+h^{\bar{p}}+h^{\bar{p}}(\k(h))^{2\bar{p}}+\big(\mu^{-1}(\k(h))\big)^{5\bar{p}-q}\bigg)
\end{align*}
and
\begin{align*}
    \EE \left( \sup_{0\leq t\leq T}|Y(t)-\bar{X}(t)|^{\bar{p}} \right) \leq C\bigg(h^{\frac{\bar{p}}{4}}+h^{\bar{p}}+h^{\bar{p}}(\k(h))^{2\bar{p}}+\big(\mu^{-1}(\k(h))\big)^{5\bar{p}-q}\bigg).
\end{align*}
In addition, it is not hard to see that
\begin{align*}
    \sup_{0\leq t\leq 1}\sup_{|x|\leq u}(|f(t,x)|\vee |g(t,x)|\vee |Lg(t,x)|)\leq 2u^5,\quad  \forall u\geq 1.
\end{align*}
So we set $\m (u)=2u^5$ and $\k(h)=h^{-\varepsilon}$, for any $\varepsilon \in (0,1/4]$. As a result, $\m^{-1}(u)=\left(u/2\right)^{1/5}$ and $\m^{-1}(\k(h))=\left(h^{-\varepsilon}/2\right)^{1/5}$. Now, choosing  $\varepsilon$ sufficiently small, choosing $p$ sufficiently large, we can derive from Theorem \ref{theorem3-2} that
\begin{align*}
    \EE \left( \sup_{0\leq t\leq 1}|Y(t)-X(t)|^{\bar{p}} \right) \leq Ch^{\bar{p}/4}
\end{align*}
and
\begin{align*}
    \EE \left( \sup_{0\leq t\leq 1}|Y(t)-\bar{X}(t)|^{\bar{p}} \right) \leq Ch^{\bar{p}/4}.
\end{align*}
which imply that the convergence order of truncated Milstein method for the time-change SDE \eqref{ex1} is $0.25$.
\par
Let us compute the approximation of the mean square error. We run M=100 independent trajectories using \eqref{numerical} for every different step sizes$10^{-1}$ $10^{-2}$, $10^{-3}$, $10^{-4}$, $10^{-5}$. We pick up $\varepsilon=0.02$, because it is hard to find the true solution for the SDE, the numerical solution with the step size $10^{-5}$ is regarded as the exact solution.
\begin{figure}[H]
    \centering
    \includegraphics[width=0.70\textwidth]{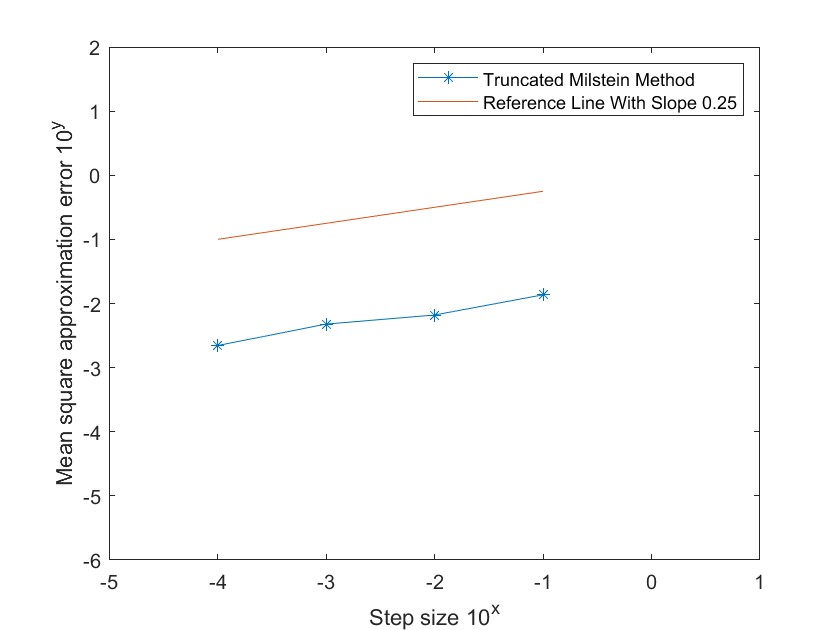}
    \caption{Convergence order of Example \ref{ex}}
\end{figure}
It is not hard to see from Figure 1 that the strong convergence order is approximately 0.25. To see it more clearly, applying the linear regression, the slope of errors against the step is 0.2517, which is quite close to the theoretical result.

\begin{expl}\label{exp}
    Consider a two-dimensional time-changed SDE
    \begin{align*}
        \label{ex2}
        \left\{
        \begin{array}{lr}
            dx_1(t)=\left([t(1-t)]^{\frac{1}{5}}x_1(t)-x_2^5(t)\right)dE(t)+\left([t(1-t)]^{\frac{1}{2}}x_2^2(t)\right)dW(E(t)),&\\
            dx_2(t)=\left([t(1-t)]^{\frac{1}{5}}x_2(t)-x_1^5(t)\right)dE(t)+\left([t(1-t)]^{\frac{1}{1}}x_1^2(t)\right)dW(E(t)).&
        \end{array}
        \right.
    \end{align*}
\end{expl}
It is clear that
\begin{align*}
    f(t,x)=
    \begin{pmatrix}
        [t(1-t)]^{\frac{1}{5}}x_1-x_2^5\\
        [t(1-t)]^{\frac{1}{5}}x_2-x_1^5
    \end{pmatrix}
    ~~~~~\text{and}~~~~~
    g(t,x)=
    \begin{pmatrix}
        [t(1-t)]^{\frac{1}{2}}x_2^2\\
        [t(1-t)]^{\frac{1}{2}}x_1^2
    \end{pmatrix}.
\end{align*}

Similar to Example 4.1, it is not hard to verify that coefficients $f(t,x)$ and $g(t,x)$ satisfy Assumption \ref{ass1} and \ref{ass5} with $\a=4$.

For any $x,y\in\RR$, it is easy to show that
\begin{align*}
        &(x-y)^{\mathrm{T}}(f(t,x)-f(t,y))+(5p-1)|g(t,x)-g(t,y)|^2\\
        = & (x_1-y_1)\bigg([t(1-t)]^{\frac{1}{5}}(x_1-y_1)-(x_2^5-y_2^{5})\bigg)+(x_2-y_2)\bigg([t(1-t)]^{\frac{1}{5}}(x_2-y_2)\\
        \quad &-(x_1^5-y_1^{5})\bigg)+(5p-1)\bigg([t(1-t)]^{\frac{1}{2}}(x_2^2-y_2^2)^2
        +[t(1-t)]^{\frac{1}{2}}(x_1^2-y_1^2)\bigg)^2\\
        \leq & (x_1-y_1)^2\left([t(1-t)]^{\frac{1}{5}}-(x_2^4+x_2^3y_2+x_2^2y_2^2+x_2y_2^3+y_2^4)\right)\\
        \quad &+(x_2-y_2)^2\left([t(1-t)]^{\frac{1}{5}}-(x_1^4+x_1^3y_1+x_1^2y_1^2+x_1y_1^3+y_1^4)\right)\\
        \quad &+2(5p-1)\bigg([t(1-t)](x_2^2-y_2^2)^2+[t(1-t)](x_1^2-y_1^2)^2\bigg)\\
        \leq & (x_1-y_1)^2\left([t(1-t)]^{\frac{1}{5}}-(x_2^4+x_2^3y_2+x_2^2y_2^2+x_2y_2^3+y_2^4)\right)\\
        \quad &+(x_2-y_2)^2\bigg([t(1-t)]^{\frac{1}{5}}-(x_1^4+x_1^3y_1+x_1^2y_1^2+x_1y_1^3+y_1^4)\bigg)\\
        \quad &+2(x_2-y_2)^{2}(5p-1)\bigg([t(1-t)](x_2+y_2)^2\bigg)+2(x_1-y_1)^{2}(5p-1)\\
        & \times\bigg([t(1-t)](x_1+y_1)^2\bigg).
\end{align*}
But
\begin{align*}
    -(x^3y+xy^3)=-xy(x^2+y^2)\leq 0.5(x^2+y^2)^2=0.5(x^4+y^4)+x^2y^2.
\end{align*}
Therefore, for any $t\in[0,1]$
\begin{align*}
        &(x-y)^{\mathrm{T}}(f(t,x)-f(t,y))+(5p-1)|g(t,x)-g(t,y)|^2\\
        \leq &(x_1-y_1)^2\bigg([t(1-t)]^{\frac{1}{5}}-0.5(x_2^4+y_2^4)+2(5p-1)[t(1-t)](x_1+y_1)^{2}\bigg)\\
        &+(x_2-y_2)^2\bigg([t(1-t)]^{\frac{1}{5}}-0.5(x_1^4+y_1^4)+2(5p-1)[t(1-t)](x_2+y_2)^{2}\bigg)\\
        \leq &C(x-y)^2,
\end{align*}
where the basic inequality $(a+b)^{2}\leq 2(a^{2}+b^{2})$ is used, and the fact that polynomials with the negative coefficients for the highest order term can always be bounds. This indicates that Assumption \ref{ass2} holds.
\par
For that Assumption \ref{ass3}, for any $q>2$ and any $t\in[0,1]$, we can drived is satisfied next
\begin{align*}
        & x^{\mathrm{T}}f(t,x)+(5q-1)|g(t,x)|^2\\
        =& ([t(1-t)]^{\frac{1}{5}}x_1^2-x_1x_2^5)+([t(1-t)]^{\frac{1}{5}}x_2^2-2x_1^5x_2)+2(5q-1)|[t(1-t)](x_1^{2}+x_2^{2})|^2\\
        \leq & [t(1-t)]^{\frac{1}{5}}(x_1^2+x_2^2)-x_1x_2(x_1^4+x_2^4)+2(5q-1)[t(1-t)](x_1^2+x_2^2)\\
        \leq & C(1+|x|^2),
\end{align*}
Then, we deal with Assumption \ref{ass4} by assuming that $\gamma_f\in (0,1]$ and $\gamma_g\in (0,1]$, for any $s,t\in[0,T]$, using the mean value theorem for the temporal variable,
\begin{align*}
        &|f(s,x)-f(t,x)|\\
        \leq& |([(s(1-s)]^{\frac{1}{5}}-[t(1-t)]^{\frac{1}{5}})x_1+([s(1-s)]^{\frac{1}{5}}-[t(1-t)]^{\frac{1}{5}})x_2|\\
        \leq& C_1|s-t|^{\frac{1}{5}}x_1+C_2|s-t|^{\frac{1}{5}}x_2,
\end{align*}
and
\begin{align*}
        &|g(s,x)-g(t,x)|\\
        \leq& |([s(1-s)]^{\frac{1}{2}}-[t(1-t)]^{\frac{1}{2}})x_2^2+([s(1-s)]^{\frac{1}{2}}-[t(1-t)]^{\frac{1}{2}})x_1^2|\\
        \leq& C_1|s-t|^{\frac{1}{2}}x_2^2+C_2|s-t|^{\frac{1}{2}}x_1^2.
\end{align*}
Thus, Assumptions \ref{ass4} is satisfied with $\gamma_f=1/5$ and $\gamma_g=1/2$. According to Theorem \ref{theorem3-2} and Example \ref{ex}, we can also set $\m (u)=2u^5$ and $\k(h)=h^{-\varepsilon}$, for any $\varepsilon \in (0,1/4]$, choosing  $\varepsilon$ sufficiently small and $p$ sufficiently large, we can derive from Theorem \ref{theorem3-1} that
\begin{align*}
    \EE \left( \sup_{0\leq t\leq 1}|Y(t)-X(t)|^{\bar{p}} \right) \leq Ch^{\bar{p}/5}
\end{align*}
and
\begin{align*}
    \EE \left( \sup_{0\leq t\leq 1}|Y(t)-\bar{X}(t)|^{\bar{p}} \right) \leq Ch^{\bar{p}/5}.
\end{align*}

which imply that the convergence order of truncated milstein method for the time-change SDE \eqref{ex1} is $1/5$ similarly. Next, we will verify through computer simulation.

Same example \ref{ex}, we run M=100 independent trajectories using \eqref{numerical} for every different step sizes $10^{-1}$,  $10^{-2}$, $10^{-3}$, $10^{-4}$, $10^{-5}$, the numerical solution with the step size $10^{-5}$ is regarded as the exact solution.

\begin{figure}[H]
    \centering
    \includegraphics[width=0.70\textwidth]{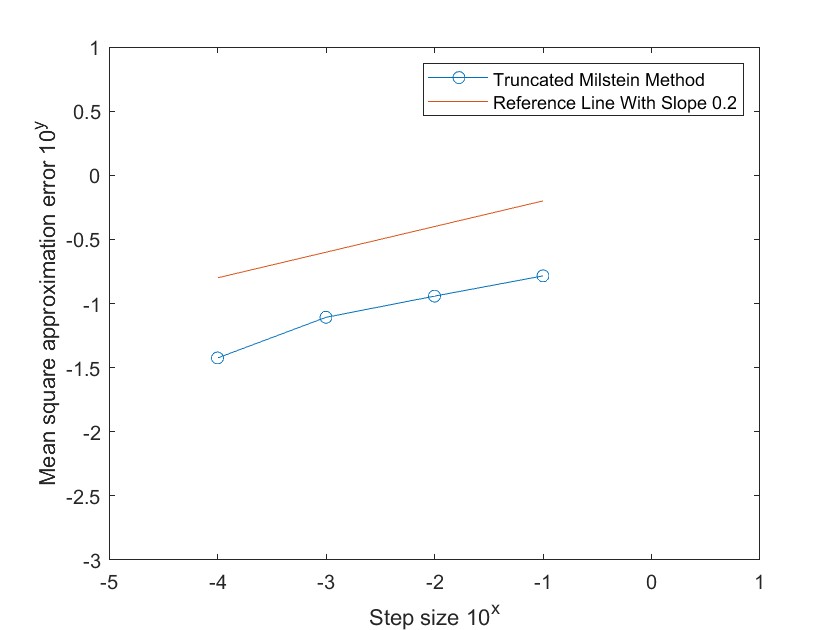}
    \caption{Convergence order of Example \ref{exp}}
\end{figure}

It is not hard to see from Figure 2 that the order of convergence can be obtained as 0.2 approximately. To see it more clearly, applying the linear regression shows that the slope of the line of errors is about 0.2086, which is also very close to the theoretical result.

\bibliography{LiuWuref}

\end{document}